\documentclass[leqno,12pt]{article}
\usepackage{latexsym}
\usepackage{amsmath}
\usepackage{amssymb}
\usepackage{bbm}
\usepackage{color}
\usepackage{graphicx}
\usepackage{epsfig}
\usepackage[latin1]{inputenc}

 0

\newtheorem{thm}{Theorem}[section]
\newtheorem{prop}[thm]{Proposition}
\newtheorem{cor}[thm]{Corollary}
\newtheorem{lem}[thm]{Lemma}
\newtheorem{defi}[thm]{Definition}
\newtheorem{rem}[thm]{Remark}

\newcommand{\ea}{\end{array}}
\newcommand{\beqohne}{\begin{eqnarray*}}
\newcommand{\eeqohne}{\end{eqnarray*}}
\newcommand{\beohne}{\begin{equation*}}
\newcommand{\eeohne}{\end{equation*}}
\newcommand{\R}{\mathbb{R}}
\newcommand{\N}{\mathbb{N}}

\newcommand{\E}{\mathbb{E}}
\def\proof{\noindent{\bf Proof:}\hskip10pt}

\def \sur#1#2{\mathrel{\mathop{\kern 0pt#1}\limits^{#2}}}

\newcommand{\la}{\lambda}

\newcommand{\muun}{\mu^{(N)}_{\u}}

\def \u{{\tt u}}

\def\proof{\paragraph{Proof:}}

\newcommand{\tr}{\mathrm{tr}}

\newcommand{\PNV}{\mathbb P_N^V}

\newcommand{\muunIj}{\mu^{(N)}_{\u,{I(j)}}}

\newcommand{\Ir}{\mathcal{J}}
\newcommand{\Fr}{\mathcal{F}}
\newcommand{\Sr}{\mathcal{S}}
\newcommand{\ap}{\alpha^+}
\newcommand{\am}{\alpha^-}
\newcommand{\bp}{b^+}
\newcommand{\bm}{b^-}

\newcommand{\SitnIj}{\tilde{\Sigma}^{(N)}_{I(j)}}
\newcommand{\Sitn}{\tilde{\Sigma}^{(N)}}

\newcommand{\lap}{\la^+}

\newcommand{\lapj}{\la^+(j)}
\newcommand{\lamj}{\la^-(j)}
\newcommand{\lapjM}{\la_M^+(j)}
\newcommand{\lamjM}{\la_M^-(j)}

\newcommand{\GUE}{\operatorname{GUE}}

\newcommand{\LUE}{\operatorname{LUE}}

\newcommand{\lal}{\langle\langle}
\newcommand{\rar}{\rangle\rangle}

\makeatletter\@addtoreset{equation}{section}\makeatother

\newcommand{\diag}{\operatorname{diag}}

\def\proof{\noindent{\bf Proof:}\hskip10pt}

\title{Sum rules and large deviations for spectral matrix measures}

\author{{\small Fabrice Gamboa}\footnote{ Universit\'e Paul Sabatier, Institut de Math\'ematiques de Toulouse,  31062 Toulouse Cedex 9, France, 
gamboa@math.univ-toulouse.fr}
\and{\small Jan Nagel}\footnote{Technische Universitat M\"unchen, Fakult\"at f\"ur Mathematik, Boltzmannstr. 3, 85748 Garching, Germany,  e-mail: jan.nagel@tum.de}
\and{\small Alain Rouault}\footnote
{Laboratoire de Mathématiques de Versailles, UVSQ, CNRS, Université Paris-Saclay, 78035-Versailles Cedex France, e-mail: alain.rouault@uvsq.fr}}

\begin{document}
\maketitle
\begin{abstract}
A sum rule relative to a reference measure on $\mathbb R$ is a relationship between the reversed Kullback-Leibler divergence of a positive measure on $\mathbb{R}$ and some non-linear functional built on spectral elements related to this measure (see for example Killip  and Simon 2003). 
In this paper, using only probabilistic tools of large deviations, we extend the sum rules obtained in  Gamboa, Nagel and Rouault (2015) 
 to the case of  Hermitian matrix-valued measures. We recover the earlier result of Damanik, Killip and Simon (2010) when the reference measure is the (matrix-valued)   semicircle law and obtain a new sum rule when the reference measure is the 
 (matrix-valued) 
 Marchenko-Pastur  law.
\end{abstract}
\bigskip

{\bf Keywords :}
{Sum rules,  block Jacobi matrix,  matrix orthogonal polynomials, matrix-valued spectral measures, large deviations, random matrices.
\smallskip

\section{Introduction}

For a probability measure $\mu$ on the unit circle $\mathbb T$ with Verblunsky coefficients $(\alpha_k)_{k \geq 0}$, Verblunsky's form of  Szeg\H{o}'s theorem may be written as
\begin{equation}
\label{Szegobasic}
\mathcal K\left(\tfrac{d\theta}{2\pi} | \mu\right) = -\sum_{k\geq 0}  \log(1 - |\alpha_k|^2)\,,\end{equation}
where $\mathcal K$ is the relative entropy, defined by
\begin{equation}
\label{KL}
{\mathcal K}(\nu\ |\ \mu)= \begin{cases}  \ \displaystyle\int_{\mathbb T}\log\frac{d\nu}{d\mu}\!\ d\nu\;\;& \mbox{if}\ \nu\ \hbox{is absolutely continuous with respect to}\ \mu, \\
   \  \infty  &  \mbox{otherwise.}
\end{cases}
\end{equation} 
Identity (\ref{Szegobasic})  is a sum rule\footnote{see the preface of \cite{simon0}} connecting an entropy and a functional of the recursion coefficients, and remains one of the most important result of the theory of orthogonal polynomials on the unit circle (OPUC) (see \cite{simon0} for extensive history and bibliography).\\
The corresponding result in the theory of orthogonal polynomials on the real line (OPRL) is the Killip-Simon sum-rule (\cite{KS03}). The reference measure is the semicircle distribution
\[\operatorname{SC}(dx) := \frac{\sqrt{4-x^2}}{2\pi} 1_{[-2, 2]} (x) dx\,.\]
The right hand side involves the Jacobi coefficients and in the left hand side appears an extra term corresponding to a contribution of isolated masses of $\mu$  outside $[-2, 2]$ (bounds states).\\ 
In a previous paper (\cite{magicrules}), we gave an interpretation and a new proof of this result from a probabilistic point of view.
This approach allowed us to prove new sum rules, when the reference measure is the Marchenko-Pastur distribution
\begin{align*}
\operatorname{MP}(\tau)(dx) = \frac{\sqrt{(\tau^+ -x)(x-\tau^-)}}
{2\pi \tau x} \, \mathbbm{1}_{(\tau^-, \tau^+)} (x) dx , 
\end{align*} 
where $\tau\in(0,1]$, $\tau^\pm = (1\pm \sqrt\tau)^2$, and also the Kesten-McKay distribution
\begin{align*}
\hbox{KMK}_{u^-, u^+}(dx) = \frac{C_{u^-, u^+}}{2\pi}\frac{\sqrt{(u^+ -x)(x- u^-)}}{x(1-x)} \ \mathbbm{1}_{(u^-, u^+)}(x) dx\, ,
\end{align*}
where $C_{u^-, u^+}$ is the normalizing constant.

Besides, known extensions of  Szeg\H{o}'s theorem (\cite{delsarte}) and of the Killip-Simon sum rule  (\cite{damal}) are available in the context of matrix-valued measures and the Matrix Orthogonal Pololynomials on the Unit Circle (MOPUC) or Matrix Orthogonal polynomials on the Real Line (MOPRL).\\ 
It seems natural to  see if the probabilistic methods are robust enough to encompass this matricial framework. Actually, the answer is positive. For the MOPUC context see \cite{GRoperatorvalued} and \cite{GNROPUC}. The aim of the present article is precisely to treat the MOPRL case. This allows to give interpretation and new proof of the Damanik-Killip-Simon's sum rule and also  to state a matrix version of the  sum rule relative to the Marchenko-Pastur distribution. \\
Let us explain the main features of our method\footnote{Note that recently Breuer et al.  \cite{BSZ} posted on arxiv a paper exposing the method to a non probabilistic audience  and giving some developments.}. As mentioned by several authors (\cite{killip2007stv}, \cite{simon2}), the main characteristic of  the above sum rules both sides are nonegative (possibly infinite) functionals. We will identify them as rate function of large deviations for random measures. Roughly speaking, that means that a sequence of random measures converges to a deterministic limit exponentially fast and the probabilities of deviating from the limit is measured by the rate function. We give two different encodings of the randomization with two rate functions $I_A$ and $I_B$, for which the uniqueness of rate functions yields the equality $I_A= I_B$. \\
To be more specific, let us give some notation. For $p \geq 1$ (fixed in all the sequel), let us denote by $\mathbb{M}_{p}$ the set of all $p\times p$ matrices with complex entries and by  $\mathcal H_p \subset \mathbb{M}_{p}$ the subset of  Hermitian matrices. 

A
matrix measure $\Sigma=(\Sigma_{i,j})_{i,j}$ of size $p$ on $\mathbb{R}$ is a matrix of signed complex measures, such that $\Sigma(A)= (\Sigma_{i,j}(A))_{i,j}\in \mathcal H_p$ for any Borel set 
$A\subset \mathbb{R}$. Further, if for any $A$,  $\Sigma(A)$ is a nonnegative matrix we say that $\Sigma$ is nonnegative.  We denote by $\mathcal{M}_p(T)$ the set of $p\times p$ nonnegative matrix measures with support in $T\subset\mathbb{R}$. Further,  $\mathcal{M}_{p,1}(T)$ is the subset of $\mathcal{M}_p(T)$ of normalized measures $\Sigma$ satisfying $\Sigma(T)={\bf 1}$, where ${\bf 1}$ is the $p\times p$ identity matrix..\\
A natural example of matrix measures comes from an application of the spectral theorem.  More precisely,  a Hermitian matrix $X$ of size $N \times N$ may be written as $UDU^*$ where $D=\diag(\lambda_i)$ contains the eigenvalues of $X$ and $U$ is the matrix formed by an orthonormal basis of  eigenvectors. 
Assume that the system $(e_1, \dots, e_p)$ of the $p$ first vectors of the canonical basis of $\mathbb C^N$ is cyclic, i.e. that Span$\{A^k e_i, k \geq 0, i \leq N\} = \mathbb C^N$. Then, for $p\leq N$, there exists a unique 
spectral
 matrix measure 
$\Sigma^X\in \mathcal{M}_p(\mathbb{R})$ supported by the spectrum of $X$, such that for all $k>0$ and 
$1 \leq i,j \leq p$
\begin{align} \label{generalspectralmeasure}
 (X^k)_{i,j}
 =  \int_{\mathbb R} x^k d\Sigma_{i,j}^X (x),
\end{align}
Let us assume that all the eigenvalues of $X$ have a single multiplicity. For $j=1,\cdots, N$, let $u^{j}:=(U_{i,j})_{i=1,\cdots,p}$ be the $j^{th}$ truncated column of $U$. Then, obviously 
\begin{align}\label{ourspectralmeasure}
\Sigma_p^X(dx) = 
\sum_{j=1}^N u_ju_j^* \delta_{\lambda_j}(dx) ,
\end{align}
where $\delta_a$ is the Dirac measure in $a$. That is, $\{\lambda_1,\dots ,\lambda_N\}$ is the support of $\Sigma^X$ and $u_1u_1^*,\dots u_Nu_N^*$ are the weights. 
Furthermore, as $U$ is unitary we have $\sum_j u_j u_j^* ={\bf 1}$ so that $\Sigma^X\in\mathcal{M}_{p,1}(T)$.\\
If we assume further that $N=pn$ for some positive integer $n$ 
 then it is possible to build a block tridiagonal matrix  
\begin{equation}
\label{jacmatrix}
J_{n}
= \begin{pmatrix} 
  B_1 & A_1    &         &         \\
                A_1^* & B_2    & \ddots  &         \\
                    & \ddots & \ddots  & A_{n-1} \\
                    &        & A_{n-1}^* & B_n
\end{pmatrix} 
\end{equation} 
such that $\Sigma^X=\Sigma^{J_{n}}$.
Here, all the blocks of $J_{n}$ are elements of $\mathbb{M}_{p}$. 
The case $p=1$ is the most classical and relies on the construction of the OPRL in  $L^2(\Sigma^X)$ (see for example \cite{simon0}). The general case is more complicated and requires  technical tools from the theory of MOPRL (see \cite{simon1}). In Section \ref{suconjac}, we will recall the construction of such tridiagonal representations.\\
 As a result, we have  two encodings of $\Sigma^X$: \eqref{ourspectralmeasure} and  \eqref{jacmatrix}.

Our measures are random in the sense that we first draw a
random matrix $X$ in $\mathcal H_N$ (with $N \geq p)$ and then considering its spectral measure $\Sigma^X$.  
In a large class of random matrix models, invariance by unitary transformations is postulated, which means that $X$ is sampled in $\mathcal H_N$ from the distribution 
\begin{align} \label{matrixensemble}
Z_N^{-1} e^{-\tr N V(X)} dX ,
\end{align}
for $V$ a confining potential and $dX$ is the Lebesgue measure on $\mathcal H_N$. Its eigenvalues behave as a Coulomb gas 
 (see formula \eqref{generaldensity}) and the matrix of eigenvectors follows independently the Haar distribution in the set of $N \times N$ unitary matrices. 
The most popular models are the Gaussian Ensemble, corresponding to the potential $V(x) = x^2/2$, the Laguerre Ensemble,  with $V(x) =  \tau^{-1} x - (\tau^{-1} - 1)\log x$ $(x > 0)$ and the Jacobi Ensemble with $V(x) = -\kappa_1  \log x - \kappa_2 \log (1-x)$ $(x \in (0,1))$.\\

We will prove in this paper that, under convenient assumptions on $V$,  as $N$ goes to infinity, the random spectral matrix measure $\Sigma^X$ converges to some equilibrium matrix measure (depending on $V$) at exponential rate. Indeed, we show that this random object satisfies a large deviation principle (LDP). To be self contained, we recall in Section \ref{sudular} the definition and useful facts on LDP.
Furthermore, the  rate function of this LDP involves a matrix extension of the reversed Kullback-Leibler information with respect to the equilibrium matrix measure (see equation (\ref{kullback}) and Theorem \ref{jointLDP}) and a contribution of the outlying eigenvalues. 



Looking for the right hand side of a possible sum rule is equivalent to look for a LDP for the encoding by means of the  sequence of blocks $A_k$ and $B_k$. 
In the scalar case ($p=1$), it is well known that
the above classical ensembles  have very nice properties \cite{dumede2002}. 
For the Gaussian  Ensemble the coefficients appearing in the tridiagonal matrix are independent with simple distributions. The diagonal terms have Gaussian distribution while the subdiagonal ones have so-called $\chi$ distributions \cite{dumede2002}. 
We will give in Lemmas \ref{triG} and \ref{triL} results in the same spirit in the general block case both for the Gaussian and Laguerre ensembles.
 These properties allow to compute the rate function of the LDP by the way of the blocks involved in the Jacobi representation. Further, the uniqueness of a rate function leads to our two main Theorems \ref{sumruleHermite} and \ref{sumruleLag} that 
are  sum rules.
 Theorem \ref{sumruleHermite} has been  proved in \cite{damal} by strong analysis tools. We recover this result by using only probabilistic arguments. For $p=1$, it has been proved earlier in  \cite{KS03}. Up to our knowledge, our Theorem 
\ref{sumruleLag} is new and is the matrix extension of the one we have obtained in \cite{magicrules} for $p=1$.

We stress that one of the main differences with purely functional analysis methods is that, for us,  the SC distribution and the free semi-infinite matrix (corresponding to  $A_k \equiv {\bf 1} $ and $B_k \equiv {\bf 0}$, with ${\bf 0}$ the $p\times p$ matrix of zeros) do not have a central role. Additionally, the non-negativity of both sides of the sum rule is automatic. 

For the Jacobi ensemble, the method used in the scalar case, based on the Szeg{\H o} mapping, is not directly extendible. So, getting a sum rule needs a more careful study. To avoid long developments here, we keep this point for a forthcoming paper.

Our paper is organized as follows. In Section \ref{smopo}, we first give definitions and tools to handle MOPRL 
 and spectral measures.
 Then, we state our main results concerning 
sum rules for spectral matrix measures.   In Section \ref{srara}, we introduce random models and state Large Deviations Principles, first  
for   random spectral measures drawn by using a general potential and   then 
  for block Jacobi coefficients in the Hermite and Laguerre cases. Section 4 is devoted to  proofs of both sum rules, up to the LDP's. 
 The most technical proofs of LDP's are postponed to Sections 
\ref{Sspopo} and  \ref{sec6}. 


\section{MOPRL 
 and block Jacobi matrices} 
\label{smopo}
\subsection{Construction}\label{sec:polynomials} 
\label{suconjac}
We will need to work with polynomials with coefficients in $\mathbb{M}_{p}$. They will be {\it orthogonal} with respect to some matrix measure on $\mathbb{R}$. 
Let us give some  notation and recall  some useful facts (see  \cite{damanik2008analytic} for details). Let 
 $\Sigma\in \mathcal M_p(\mathbb{R})$ be a compactly supported matrix measure. 
Further, let  $F$ and $G$ be 
continuous matrix valued  functions $ F,G: \mathbb R \rightarrow \mathbb M_{p}$. We define the two pseudo-scalar products $p\times p$ 
(which are elements of $\mathbb M_p$)   by setting
\begin{eqnarray*}\lal F, G\rar_R &=& \int_\mathbb R F(z)^* d\Sigma(z)G(z)\,, \\
\lal F, G\rar_L &=& \int_\mathbb R G(z) d\Sigma(z)F(z)^*\,.
\end{eqnarray*}
A sequence of functions on $\mathbb{R}$,  $(\varphi_j)$ with values in $\mathbb M_p$  is called right-orthonormal if
\[\lal\varphi_i, \varphi_j\rar_R= \delta_{ij}{\bf 1} \ .\]
The orthogonal polynomial recursion is built as follows. 
First, assume that $\Sigma$ is nontrivial, that is, 
\begin{align}\label{nontrivial}
\tr \lal P,P \rar _R>0
\end{align} 
for any non zero matrix polynomial $P$ (see Lemma 2.1 of \cite{damanik2008analytic} for equivalent characterizations of nontriviality). 
Applying the 
Gram-Schmidt procedure  to $\{{\bf 1}, x{\bf 1}, x^2 {\bf 1},\dots\}$, we obtain a sequence $(P_n^R)_n$  of  right monic matrix orthogonal polynomials. In other words, $P_n^R$ is the unique matrix polynomial $ P_n^R(x) = x^n{\bf 1} +$ lower order terms such that $\lal x^k{\bf 1}, P_n^R\rar_R ={\bf 0}$ for $k=0, \dots, n-1$. 
For nontrivial matrix measures, this is possible for any $n\geq 0$ and this sequence satisfies the recurrence relation
\begin{equation}
\label{relrecR}
xP_n^R(x) = P_{n+1}^R (x) + P_n^R (x)u_n^R + P_{n-1}^R v_n^R\,. 
\end{equation}
If we set
  \begin{eqnarray*}  \gamma_n := \lal P_n^R, P_n^R\rar_R , 
\end{eqnarray*}
then $\gamma_n$ is positive definite and we have
\[v_n^R = \gamma_{n-1}^{-1} \gamma_n \,.\]
To get normalized orthogonal polynomials $p_n^R$ we set
\begin{equation}\label{smallp}p_0^R = {\bf 1}\ \  ,\ \  p_n^R = P_n^R\kappa_n^R\end{equation}
where for every $n$, $\kappa_n^R \in \mathbb M_p$ 
 has to satisfy 
\begin{equation}
\gamma_n =  \left(\kappa_n^R\left(\kappa_n^R\right)^ *\right)^{-1}\,.
\end{equation}
This constraint opens several choices for $\kappa_n^R$ (see Section 2.1.5 of \cite{damanik2008analytic}). Let us leave the choice open, setting
\begin{equation}
\label{R1}
\kappa_n^R = \gamma_n^{-1/2} \sigma_n 
\end{equation}
with $\sigma_n$ unitary and $\sigma_0 = {\bf 1}$.

\begin{rem}We can define similarly the sequence of monic polynomials $P_n^L$ and the sequence of left-orthonormal polynomials $p_n^L$ in the same way.
We have
\begin{eqnarray*}P_n^R = (P_n^L)^*\end{eqnarray*}
and  \begin{eqnarray*}   \lal P_n^L , P_n^L \rar_L  =  \gamma_n\end{eqnarray*}
and the recurrence relation:
\begin{eqnarray}
 xP_n^L (x) = P_{n+1}^L (x) +  u_n^L P_n^L(x) + v_n^L P_{n-1}^L
\end{eqnarray}
with
\[ v_n^L =  \gamma_n  \gamma_{n-1}^{-1}\,.\]
The above condition (\ref{R1}) is replaced by  
$p_n^L = \tau_n \gamma_n^{-1/2} P_n^L$.
%
\end{rem}

To formulate the recursion in terms of orthonormal polynomials, we use (\ref{smallp}) and get
\begin{eqnarray}xp_n^R &=& p_{n+1}^R ( \kappa_{n+1}^R)^{-1}\kappa_n^R + p_n^R(\kappa_n^R)^{-1}u_n^R \kappa_n^R + p_{n-1}^R(\kappa_{n-1}^R)^{-1}v_n^R \kappa_n^R\end{eqnarray}
i.e.
\begin{eqnarray}\label{polrek}
xp_n^R &=& p_{n+1}^R A_{n+1}^* + p_n^R B_{n+1} + p_{n-1}^R A_n\end{eqnarray}
with
\begin{eqnarray}
\label{a1}
A_n &=& (\kappa_{n-1}^R)^{-1}v_n^R \kappa_n^R = \sigma_{n-1}^* \gamma_{n-1}^{1/2}v_n^R \gamma_n^{-1/2}\sigma_n = \sigma_{n-1}^* \gamma_{n-1}^{-1/2}\gamma_n^{1/2}\sigma_n\\
\nonumber
 B_{n+1} &=& (\kappa_n^R)^{-1}u_n^R \kappa_n^R = \sigma_n^* \gamma_n^{1/2}u_n^R \gamma_n^{-1/2}\sigma_n\,.
\end{eqnarray}
Note that \eqref{a1} yields
\begin{equation}
\label{asquare}
A_nA_n^* = \sigma_{n-1}^* \gamma_{n-1}^{-1/2}\gamma_n \gamma_{n-1}^{-1/2}\sigma_{n-1}\,.
\end{equation}

In other terms the map $f \mapsto (x \mapsto xf(x))$ defined on the space of matrix polynomials is a right homomorphism and is represented in the (right-module) basis $\{p_0^R,p_1^R,\dots \}$ by the matrix
\begin{align} \label{jacobimatrix}
J = 
 \begin{pmatrix} 
  B_1 & A_1    &                  \\
                A_1^* & B_2    & \ddots           \\
                    & \ddots & \ddots   
\end{pmatrix}
\end{align}
with $B_k$ Hermitian and $A_k$ non-singular. Moreover 
 the measure $\Sigma$ is again the spectral 
 measure of the matrix $J$ defined as in \eqref{generalspectralmeasure} (Theorem 2.11 of \cite{damanik2008analytic}). Let us remark that although to each $\Sigma$ corresponds a whole equivalence class of Jacobi coefficients given by the different $\sigma_n$, there is exactly one representative such that all $A_k$ are Hermitian positive definite 
(Theorem 2.8 in \cite{damanik2008analytic}).

Starting with a finite dimensional Jacobi matrix $J_n$ as in \eqref{jacmatrix}, the spectral matrix measure of $J_n$ is supported by at most 
$n$
points and is in particular not nontrivial. However, we may still define $p_1^R,\dots ,p_{n-1}^R$ by the recursion \eqref{polrek}. As long as the 
$A_k$'s are invertible, these polynomials are orthonormal with respect to the spectral measure of $J_n$ and \eqref{nontrivial} holds for all 
 polynomials up to degree $n-1$. 


\medskip

If $\Sigma$ is a {\it quasi scalar} measure, that is  if $\Sigma =  \sigma \cdot {\bf 1}$ with $\sigma$ a scalar measure and if $\Pi$ is a positive matrix measure with Lebesgue decomposition
\[\Pi(dx) = h(x) \sigma(dx) +  \Pi^s(dx)\,,\] we extend the definition  (\ref{KL}) by
\begin{equation}
\label{kullback}
\mathcal K(\Sigma | \Pi) := - \int \log\det h(x)\ \sigma (dx)\,.
\end{equation}
\begin{rem}
It is possible to rewrite the above quantity in the flavour of Kullback-Leibler information (or relative entropy) with the notation of \cite{mandrekar1} or \cite{robertson}, i.e.
\[\mathcal K(\Sigma | \Pi) = \int \log \det \frac {d\Sigma(x)}{d\Pi (x)} d\sigma(x)\,,\]
if $\Sigma$ is {\it strongly absolutely continuous} with respect to $\Pi$, and infinity otherwise. 
See Corollary 8 in \cite{GRoperatorvalued}.
\end{rem}

\subsection{Measures on $[0, \infty)$}
\label{sec:polynomialson0infty}
When the measures are supported by $[0, \infty)$, there is a specific form of the Jacobi coefficients, leading to a particularly interesting parametrization, which will be crucial in the Laguerre ensemble.

In \cite{destu02}, it is proved that if a nontrivial matrix measure $\Sigma$ has a support included in $[0, \infty)$ then there exists a sequence $(\zeta_n)_n$ of non-singular elements of $\mathbb M_{p}$ such that
\begin{equation}
\label{the_zetas} u_n^R = \zeta_{2n+1} + \zeta_{2n} \ \ , \ \  v_n^R = \zeta_{2n-1}\zeta_{2n} \end{equation}
with $\zeta_0={\bf 0}$ and moreover, 
\begin{equation} \label{zetadecomposition}
\zeta_n=h_{n-1}^{-1}h_n, (n \geq 1) 
\end{equation}
with $h_0={\bf 1}$ and for $n\geq 1$,  $h_n \in \mathcal H_p$ is 
 positive definite.
 Note that this implies \[\gamma_n=\zeta_1\dots \zeta_{2n}=h_{2n}\,.\]
From (\ref{a1}) we then have the representations of the Jacobi coefficients
\begin{align} 
\label{repB}B_{n+1} & = \sigma_n^* \gamma_n^{1/2}(\zeta_{2n+1} +\zeta_{2n})\gamma_n^{-1/2}\sigma_n \\
\nonumber
A_n &= \sigma_{n-1}^*\gamma_{n-1}^{1/2}\zeta_{2n-1}\zeta_{2n}\gamma_{n}^{-1/2}\sigma_{n}
\,.\end{align}
In the scalar case, this yields $B_1 = \zeta_1$ and for $n \geq 1$
\begin{equation}
\label{scal}
B_{n+1}=  \zeta_{2n+1} + \zeta_{2n}  , \qquad (A_n)^2 = \zeta_{2n-1}\zeta_{2n}\,. 
\end{equation}
In the matrix case, 
we may set
\begin{equation}
\label{defZZ}
Z_{2n+1} = \sigma_n^*\gamma_n^{1/2} \zeta_{2n+1}\gamma_n^{-1/2}\sigma_n \ \ , \ \ Z_{2n} = \sigma_n^*\gamma_n^{1/2}\zeta_{2n}\gamma_n^{-1/2}\sigma_n\,.
\end{equation}
To highlight a further decomposition, we set 
\begin{align*}
C_n=\sigma_n^*h_{2n}^{1/2}h_{2n-1}^{-1/2},\qquad D_{n+1}=\sigma_n^* h_{2n}^{-1/2} h_{2n+1}^{1/2} \,.
\end{align*}
With these definitions and \eqref{zetadecomposition} we see that the matrices 
\begin{equation} \label{scalmat}
Z_{2n+1}=D_{n+1}D_{n+1}^*,\qquad Z_{2n}=C_nC_n^*
\end{equation}
are in fact Hermitian positive definite. For the recursion coefficients we get the following matrix analogues of \eqref{scal}, 
$B_1 = D_1D_1^*$,  and for $n \geq 1$ 
\begin{equation}
\label{ABfromZ}B_{n+1}=  D_{n+1}D_{n+1}^* + C_nC_n^* =Z_{2n+1} + Z_{2n} , \qquad A_n  = D_nC_n^*\,.
\end{equation}
In other words, the Jacobi operator $J$ can in fact be decomposed as $J=XX^*$, where $X$ is the bidiagonal matrix
\begin{equation} \label{bidiagonal}
X = \begin{pmatrix} 
  D_1 & {\bf 0}    &                 \\
                C_1 & D_2    & \ddots         \\
                    & \ddots & \ddots   \\
\end{pmatrix} .
\end{equation}
Moreover, the entries of $X$ can be chosen to be Hermitian positive definite. We are still free to choose the unitary matrices $\sigma_n$ (although we have to fix $\sigma_1={\bf 1}$) in the definition of orthonormal polynomials and we let $U$ be the block-diagonal matrix with $\sigma_1,\sigma_2,\dots$ on the diagonal. Moreover, let $P$ denote a block-diagonal matrix with unitary $p\times p$ matrices $\tau_1,\tau_2,\dots$ on the diagonal. Then our measure $\Sigma$ is also the spectral matrix measure of 
 $UXPP^*X^*U^* = (UXP)(UXP)^*$. 
The matrix $UXP$ has the form
\begin{align*}
UXP = \begin{pmatrix} 
  \sigma_1D_1\tau_1 & {\bf 0}    &              &     \\
                \sigma_2C_1\tau_1 & \sigma_2D_2\tau_2    & {\bf 0} &            \\
      &    \sigma_3C_2\tau_2            &  \sigma_3D_3\tau_3  &  \ddots &  \\
      &   &  \ddots & \ddots & \\ 
\end{pmatrix} .
\end{align*}
For the first entry, $\sigma_1={\bf 1}$ and $D_1$ is always Hermitian positive definite, so we may set $\tau_1={\bf 1}$. Recall that for $A$ a non-singular matrix, there exists a unique unitary $\sigma$ such that $A\sigma$ is Hermitian positive definite, and if $\Sigma$ is nontrivial, all $D_k,C_k$ are non-singular. Therefore, we can recursively choose $\sigma_{k+1}$ such that $\sigma_{k+1}C_k\tau_k$ is Hermitian positive definite and then $\tau_{k+1}$ such that $\sigma_{k+1}D_{k+1}\tau_{k+1}$ is positive definite. This yields a unique decomposition with positive definite blocks.

\subsection{Sum rules}

For $\alpha^-<\alpha^+$, let $\mathcal{S}_p = \mathcal{S}_p(\am,\ap)$ be the set of all bounded nonnegative 
measures $\Sigma \in\mathcal{M}_p(\mathbb{R})$ with 
\begin{itemize}
\item[(i)] $\operatorname{supp}(\Sigma) = K \cup \{\lambda_i^-\}_{i=1}^{N^-} \cup \{\lambda_i^+\}_{i=1}^{N^+}$, where $K\subset I= [\am,\ap]$, $N^-,N^+\in\N\cup\{\infty\}$ and 
\begin{align*}
\lambda_1^-<\lambda_2^-<\dots <\am \quad \text{and} \quad \lambda_1^+>\lambda_2^+>\dots >\ap .
\end{align*}
\item[(ii)] If $N^-$ (resp. $N^+$) is infinite, then $\lambda_j^-$ converges towards $\am$ (resp. $\lambda_j^+$ converges to $\ap$).
\end{itemize}
Such a measure $\Sigma$ can be written as
\begin{align}\label{muinS0}
\Sigma = \Sigma_{|I} +  \sum_{i=1}^{N^+} \Gamma_i^+ \delta_{\lambda_i^+} + \sum_{i=1}^{N^-} \Gamma_i^- \delta_{\lambda_i^-},
\end{align}
for some nonnegative Hermitian matrices $\Gamma_1^+,\cdots, \Gamma_{N^+}^+,\Gamma_1^-,\cdots, \Gamma_{N^-}^-$.
Further, we define $\mathcal{S}_{p,1}=\Sr_{p,1}(\am,\ap):=\{\Sigma \in \mathcal{S}_p(\alpha^-,\alpha^+) |\, \Sigma(\R)={\bf 1}\}$.

\subsubsection{The Hermite case revisited}

In the scalar frame ($p=1$), the Killip-Simon sum rule gives two  different expressions for the divergence between a probability measure and the semicircle distribution (see \cite{KS03} and \cite{simon2}, Chapter 3). In the more general case $p\geq 2$, it 
gives two forms for the divergence with respect to 
\begin{align*}
\Sigma_{SC} = \operatorname{SC} \cdot {\bf 1} \ \ \hbox{with} \ \operatorname{SC}(dx) = \frac{\sqrt{4-x^2}}{2\pi} \, \mathbbm{1}_{[-2, 2]}(x)\!\ dx
\end{align*}
supported by $[\alpha^-,\alpha^+]=[-2,2]$. 
We refer to \cite{damal} Formula (10.4) and \cite{simon2}, Formula (4.6.13) for this matrix sum rule.
 The 
 block Jacobi 
matrix associated with $\Sigma_{SC}$ has entries
\begin{align*}
B_k^{\operatorname{SC}} = {\bf 0},\quad A_k^{\operatorname{SC}} = {\bf 1},\quad 
\end{align*}
for all $k\geq 1$. The spectral side of the sum rule involves a contribution of outlying eigenvalues, for which we define
\begin{align*}
\mathcal{F}_H^+(x) :=  \begin{cases}   \displaystyle \int_2^x \sqrt{t^2-4}\!\ dt = \tfrac{x}{2} \sqrt{x^2-4} - 2 \log \left( \tfrac{x+\sqrt{x^2-4}}{2}\right) & \mbox{ if} \ x \geq 2, \\
 \infty & \mbox{ otherwise}
\end{cases}
\end{align*}
and $\mathcal{F}_H^-(x)=\mathcal{F}_H^+(-x)$. Let $G$ be  the very popular function (Cram\' er transform of the exponential distribution)
\[G(x) = x - 1 - \log x \ \ (x > 0)\,.\]
We adopt the convention of the functional calculus, so that  for $X\in \mathcal H_p$ positive, 
 we have
\begin{equation}
\label{functcalc}
\tr\!\ G(X)= \tr X-\log \det X-p\,.
\end{equation}
Here is the first sum rule. 
 This remarkable equality has been first proven by \cite{damal}. In Section \ref{sec:LDP2SR}, we give a probabilistic proof. Indeed,  we show that this sum rule is a consequence of two large deviation results. 

\medskip

\begin{thm}\label{sumruleHermite}
Let $\Sigma \in \mathcal{M}_{p,1}(\R)$ be a spectral measure with Jacobi 
matrix  \eqref{jacobimatrix}. If $\Sigma \in \mathcal{S}_{p,1}(-2,2)$, then 
\begin{align*}
\mathcal{K}(\Sigma_{SC}\!\ |\!\ \Sigma) + 
\sum_{k=1}^{N^+} {\mathcal F}_{H}^+ (\lambda_k^+)  +  \sum_{k=1}^{N^-} {\mathcal F}_{H}^- (\lambda_k^-)
= \sum_{k=1}^\infty  \tfrac{1}{2} \tr  B_k^2  + \tr\, G(A_kA_k^*),
\end{align*}
where both sides may be infinite simultaneously. If $\Sigma \notin \mathcal{S}_{p,1}(-2,2)$, the right hand side equals $+\infty$.
\end{thm}

\medskip

We remark that since $\tr\, G(\sigma AA^*\sigma^*) = \tr\, G(AA^*)$ for any unitary $\sigma$, the value of the right hand side in Theorem \ref{sumruleHermite} is independent of the choice of 
$\sigma_n$'s in \eqref{a1}.

Let us restate the sum rule as in the notation of \cite{simon2}. 
To a spectral measure $\Sigma$ supported by $[-2,2]$, we associate the $m$-function (Stieltjès transform)
\begin{align*}
m(z) = \int \frac{1}{x-z} d\Sigma(x),\; (z\in\mathbb{C}\setminus[-2,2]).
\end{align*}
For $z\in \mathbb D$ (interior of the unit disk), the function $M(z) = -m(z +z^{-1})$ admits radial limits : 
 for almost all $\theta \in [0,2\pi]$, the limit 
$M(e^{i\theta}) = \lim_{r\uparrow 1} M(re^{i\theta})$ exists and is neither vanishing nor infinite. Finally, let
\begin{equation}
\label{defQ}
 {\mathcal Q} (\Sigma) = \frac{1}{\pi} \int_0^{2\pi} \log \left( \frac{\sin^p \theta }
{\det (\operatorname{Im} M(e^{i\theta}))} \right) \sin^2 \theta\!\  d\theta \,.
\end{equation}
Then the following statement is the combination of Theorem 4.6.1 and Theorem 4.6.3 in Simon's book  \cite{simon2}. It is a {\it gem}, as defined on p.19 of \cite{simon2}.
 
\medskip

\begin{thm} \label{sumrule}
Let $A_n, B_n$ be the 
 entries of the block Jacobi 
matrix $J$ and let $\Sigma$ denote the spectral 
 measure of $J$. Then
\begin{align*}
\sum_{n=1}^\infty \tr ((A_nA_n^*)^{1/2} -I)^2 + \tr B_n^2 < \infty
\end{align*}
if and only if 
\begin{itemize}
\item[(a)] The essential support of $J$ satisfies
\begin{align*} \sigma_{ess}(J) \subset [-2,2] 
\end{align*}
\item[(b)] The eigenvalues $\{ \lambda_j \}_{j=1}^\infty \notin \sigma_{ess}(J)$ satisfy 
\begin{align*} \sum_{k=1}^\infty (|\lambda_k| -2)^{3/2} < \infty 
\end{align*}
\item[(c)] 
If $\Sigma$ admits the decomposition
\[d\Sigma(x) = f(x) dx + d\Sigma_s (x)\]
with $f \in$ and $d\Sigma_s$ singular with respect to $dx$, then
\begin{align*} \int_{-2}^2 \sqrt{4-x^2} \log \det f(x) dx > -\infty\,. 
\end{align*}
\end{itemize}
In this case, we have 
\begin{equation}
\label{newsum1}
\sum_{n=1}^\infty \left( \tfrac{1}{2} \tr (B_n^2) + \tr(G(A_nA_n^*)) \right) =  {\mathcal Q} (\Sigma) + \sum_{\lambda \notin \sigma_{ess}(J)}  F(\lambda) .
\end{equation}
Here, $F$ is defined to be equal to $\mathcal{F}_H^+$ on $[2,\infty)$ and equal to $\mathcal{F}_H^-$ on $(-\infty,2]$.
\end{thm}

\medskip

Note that the integral ${\mathcal Q}(\Sigma)$ appearing in \eqref{newsum1} and defined in (\ref{defQ}) can be interpreted as a relative entropy, like in the scalar case treated in \cite{simon2}, Lemma 3.5.1. For a measure $\Sigma$ supported on $[-2,2]$, the (inverse) Szeg\H{o}\ mapping pushes forward $\Sigma$ to a measure $\Sigma_{Sz}$ on the unit circle symmetric 
with respect to complex conjugation, such that for all measureable and bounded $\varphi$,
\begin{align*}
\int_0^\pi \varphi (2\cos(\theta)) d\Sigma_{Sz}(e^{i\theta}) = \int_{-2}^2 \varphi\left( x \right) d\Sigma(x) .
\end{align*}
A straightforward generalization of the arguments in the above reference 
 show
\begin{align*}
 \operatorname{Im} M(e^{i\theta}) =  \operatorname{Im} m(2\cos \theta) = \pi f(2\cos \theta)
\end{align*}
for $\theta \in [0,\pi]$.
 Then, using the symmetry 
 and setting $x=2\cos \theta$, we obtain
\begin{align*}
   {\mathcal Q}(\Sigma) =& \frac{2}{\pi} \int_{-2}^2 \log \left( \frac{2^{-p} (4-x^2)^{-p/2}}{\det \pi f(x)}\right) \frac{1}{4} \sqrt{4-x^2}\!\ dx \\
=&  \int_{-2}^2 \log \det \left( \frac{1}{2\pi} \sqrt{4-x^2} f(x)^{-1} \right) \frac{1}{2\pi} \sqrt{4-x^2}\!\ dx \\
=&  \int \log \det \left( \frac{d\Sigma_{SC}}{d\Sigma} \right) d\operatorname{SC}
= \mathcal K(\Sigma_{SC}\!\ |\!\ \Sigma) .
\end{align*}

\medskip

\subsubsection{Our new sum rule: the Laguerre case}

In the Laguerre case, the central measure is the matrix Marchenko-Pastur law with scalar version
\begin{align*}
\operatorname{MP}(\tau)(dx) = \frac{\sqrt{(\tau^+ -x)(x-\tau^-)}}
{2\pi \tau x} \, \mathbbm{1}_{(\tau^-, \tau^+)} (x) \!\ dx , 
\end{align*} 
where $\tau\in(0,1]$, $\alpha^\pm=\tau^\pm = (1\pm \sqrt\tau)^2$ and we set $\Sigma_{\operatorname{MP}(\tau)} = \operatorname{MP}(\tau) \cdot {\bf 1}$.  The block Jacobi matrix associated with $\Sigma_{\operatorname{MP}}(\tau)$ has entries:
\[A_k^{MP} = \sqrt{\tau}\cdot {\mathbf 1}\ ,  \ (k \geq 1)\,, \ B_1^{MP} = {\mathbf 1} \ , \ B_k^{MP}= (1 +\tau)\cdot {\mathbf 1} \ , \ (k\geq 2)\,,\]
which corresponds to $\zeta_{2k-1} = {\mathbf 1} \ , \ \zeta_{2k} = \tau\cdot {\mathbf 1}$. 

For the new Laguerre sum rule, we have to replace $\mathcal{F}_H^\pm$ by
\begin{align*}
\mathcal{F}_L^+(x) =  \begin{cases}  \displaystyle \int_{\tau^+}^x \frac{\sqrt{(t -  \tau^-)(t - \tau^+)}}{t\tau}\!\ dt  & \mbox{ if} \ x \geq \tau^+,\\
\infty & \mbox{ otherwise}
\end{cases}
\end{align*}
and 
\begin{align*} 
{\mathcal F}_L^-(x) =  \begin{cases}  \displaystyle\int_x^{\tau^-} \frac{\sqrt{(\tau^- -t)(\tau^+ -t)}}{t\tau}\!\ dt  & \mbox{if} \ x \leq \tau^-,\\
      \infty & \mbox{ otherwise.}
\end{cases}  
\end{align*}
One of our main results is Theorem \ref{sumruleLag}. Up to our knowledge, this result is new.
The proof is again in Section \ref{sec:LDP2SR}. 

\medskip

\begin{thm} \label{sumruleLag}
Assume the Jacobi matrix $J$ is nonnegative definite and let $\Sigma$ be the spectral 
 measure associated with $J$. Then for any $\tau \in (0,1]$, if $\Sigma \in \mathcal{S}_{p,1}(\tau^-,\tau^+)$, 
\begin{align*} 
 {\mathcal K}(\Sigma_{MP(\tau)} 
\!\ |\!\ \Sigma) +  \sum_{n=1}^{N^+} {\mathcal F}_L^+(\lambda_n^+)  +  \sum_{n=1}^{N^-} {\mathcal F}_L^-(\lambda_n^-)
= \sum_{k=1}^\infty \tau^{-1}\tr\!\ G(\zeta_{2k-1}) +  \tr\!\ G(\tau^{-1}\zeta_{2k})
\end{align*}
where both sides may be infinite simultaneously and $\zeta_k$ is defined as in \eqref{the_zetas}. If $\Sigma \notin \Sr_{p,1}(\tau^-,\tau^+)$, the right hand side equals $+\infty$. 
\end{thm}

\medskip

Since the matrices $\zeta_k$ can be decomposed as in \eqref{zetadecomposition}, they are in fact similar to a Hermitian matrix 
$$
\zeta_k = h_{k-1}^{-1/2}\left[h_{k-1}^{-1/2} h_k h_{k-1}^{-1/2}\right]h_{k-1}^{1/2}\,,
$$ hence the sum on the right hand side in Theorem \ref{sumruleLag} is real valued.

Similar to the matrix {\it gem} , Theorem \ref{sumrule}, we can formulate equivalent conditions on the matrices $\zeta_k$ and the spectral measure, which characterize finiteness of the sum. The following corollary is the matrix counterpart of Corollary 2.4 in \cite{magicrules}. It follows immediately from Theorem \ref{sumruleLag}, since $\mathcal F_L^- (0) = \infty$ and 
\[\mathcal F^\pm_L (\tau^\pm \pm h) = \frac{4}{3\tau^{3/4}(1 \pm\sqrt \tau)^2 }h^{3/2} + o(h^{3/2}) \ \ \ (h \rightarrow 0^+)\]
and, for $H$ similar to a Hermitian matrix,
\[\tr\, G(
{\bf 1}
+H) = \frac{1}{2} \tr H^2 + o(||H||) \ \ \ \ (||H|| \rightarrow 0).\]
Here, $||\cdot||$ is any matrix norm.

\begin{cor} \label{semigemL}
Assume the Jacobi matrix $J$ is nonnegative definite and let $\Sigma$ be the spectral 
 measure of $J$. Then
\begin{align}
\label{zl2}
\sum_{k=1}^\infty[\tr(\zeta_{2k-1}-{\bf 1})^2 + \tr(\zeta_{2k} - \tau {\bf 1})^2] < \infty
\end{align}
if and only if
\begin{enumerate}
\item 
  $\Sigma \in \mathcal S_{p,1} (\tau^-, \tau^+)$
\item $\sum_{i=1}^{N^+} (\lambda_i^+ - \tau^+)^{3/2} + \sum_{i=1}^{N^-} (\tau^- - \lambda_i^- )^{3/2}  < \infty$ and if $N^->0$, then $\lambda_1^- > 0$. 

\item the spectral 
 measure $\Sigma$ of $J$ with Lebesgue decomposition $d\Sigma(x) = f(x) dx+d\Sigma_s(x)$ 
satisfies
\begin{align*}
\int_{\tau^-}^{\tau^+} \frac{\sqrt{(\tau^+-x)(x-\tau^-)}}{x} \log \det (f(x)) dx >-\infty .
\end{align*}
\end{enumerate}
\end{cor}


\section{Randomization and large deviations}
\label{srara}
\subsection{Matrix random models}

The results of the previous section rely on
two classical distributions of random Hermitian matrices: the Gaussian (or Hermite) and the Laguerre (or Wishart) ensemble. 
We denote by $\mathcal N(0, \sigma^2)$ the centered Gaussian distribution with variance $\sigma^2$. 
A random variable $X$ taking values in $\mathcal H_N$ is distributed according to the Gaussian unitary ensemble $\GUE_N$, if all real diagonal entries are 
distributed as $\mathcal (0, 1)$ 
 and the if the real and imaginary parts of off-diagonal variables are independent and $\mathcal N(0, 1/2)$ distributed (this is called complex standard normal distribution). 
All entries are assumed to be independent up to symmetry and conjugation.
The random matrix $X/\sqrt N$ has then the distribution given by (\ref{matrixensemble}) and the joint density of the (real) eigenvalues $\lambda = (\lambda_1,\dots ,\lambda_N)$ of $X$
 is 
 (see for example \cite{agz})
\begin{align} \label{evg}
g_G(\lambda) = c_G \prod_{1\leq  i < j\leq N} |\lambda_i - \lambda_j|^2  \prod_{i=1}^N e^{- \lambda_i^2/2}.
\end{align}
In analogy to the scalar $\chi^2$ distribution, the Laguerre ensemble is the distribution of the square of Gaussian matrices as follows.  If $G$ denotes a $N \times \gamma$ matrix with independent complex standard normal entries, then $GG^*$
is said to be distributed according to the Laguerre ensemble $ \LUE_N(\gamma)$ with parameter $\gamma$. If $\gamma \geq N$, the eigenvalues have the density 
(see for example \cite{agz})
\begin{align} \label{evl}
g_L(\lambda) = c_L^\gamma \prod_{1\leq  i < j\leq N} |\lambda_i - \lambda_j|^2  \prod_{i=1}^N \lambda_i^{\gamma-N} e^{-\lambda_i} \mathbbm{1}_{\{ \lambda_i>0\} }.
\end{align}

It is a well-known consequence of the invariance under unitary conjugation, that in the classical ensembles \eqref{matrixensemble}, the array of random eigenvalues and the random eigenvector (unitary) matrix are independent. Further, this  latter matrix is Haar distributed (\cite{dawid1977}). This implies the following equality in distribution for the weights  given in Lemma \ref{weightdistr} (see Proposition 3.1 in \cite{GRoperatorvalued}), which is a matrix version of the beta-gamma relation for scalar random variables. First we need a definition.

\begin{defi}
\label{betagamma}
\begin{enumerate}
\item
If $v_1,\dots ,v_N$ are independent complex standard normal distributed vectors in $\mathbb{C}^p$, set $V_j = v_jv_j^*$  for $j \leq N$. 
  We say that 
$(V_1,\dots ,  V_N)$ follows the distribution $\mathbb G_{p,N}$ on $(\mathcal H_p)^N$. 

\item If $U$ be Haar distributed in the set of $N\times N$ unitary matrices, set   $u_j = (U_{i,j})_{1\leq i\leq p}\in \mathbb{C}^p$ and $W_j = u_ju_j^*$ for $j \leq N$. 
We say that  $(W_1, \dots, W_N)$ 
 follows the  $\mathbb D_{p,N}$ distribution  on $(\mathcal H_p)^N$.
\end{enumerate}
\end{defi}

In the scalar case, the array of weights $W_j= u_j^2$ is uniformly distributed on the simplex $\{w_1 , \dots w_N \in [0,1] : \sum_i w_i = 1\}$.
\medskip

\begin{lem} \label{weightdistr}
If $(V_1,\dots ,  V_N)$ follows the distribution $\mathbb G_{p,N}$ and if 
 $H=\sum_{k=1}^N V_k$,
then 
$\big( H^{-1/2} V_1H^{-1/2},\dots , H^{-1/2} V_N H^{-1/2}\big) $ follows the distribution $\mathbb D_{p,N}$.
\end{lem}

\medskip

Our first large deviation principle will hold for a general class of $p\times p$ matrix measures. We draw the random eigenvalues $\lambda_1,\dots ,\lambda_N$ from the absolute continuous distribution $\PNV$
\begin{align}\label{generaldensity}
d\PNV (\lambda) = \frac{1}{Z_V^N} \prod_{1\leq  i < j\leq N} |\lambda_i - \lambda_j|^2
\prod_{i=1}^N e^{- N V(\lambda_i)} .
\end{align}
We suppose that the potential $V$ is continuous and real valued on the interval ${(b^-, b^+)}$ ($-\infty\leq b^-<b^+\leq+\infty$), infinite outside of $[b^-,b^+]$ and $\lim_{x\to b^\pm} V(x) = V(b^\pm)$ with possible limit $V(b^\pm)=+\infty$. Under the assumption 
\begin{itemize}
\item[(A1)] Confinement: 
$\displaystyle \qquad
\liminf_{x \rightarrow b^\pm} \frac{V(x)}{ \log |x|} > 2 \, , $
\end{itemize}
the empirical distribution $\muun$ of eigenvalues $\lambda_1,\dots ,\lambda_N$ has a 
 limit $\mu_V$ in probability, which is the unique minimizer of 
\begin{align}
\label{ratemuu}
\mu  \mapsto \mathcal E (\mu) := \int V(x) d\mu(x) - \iint \log |x-y| d\mu(x)d\mu(y).
\end{align}
and which has a compact support (see \cite{johansson1998fluctuations} or \cite{agz}). 
This convergence  can be viewed as a consequence of the LDP for the sequence $(\mu_\u^{(N)})_N$. 
We need two additional assumptions on $\mu_V$:
\begin{itemize}
\item[(A2)] One-cut regime: the support of $\mu_V$ is a single interval $[\alpha^-, \alpha^+]\subset [b^-, b^+]$ ( $\alpha^-< \alpha^+$).
\item[(A3)] Control (of large deviations): the effective potential
\begin{align}
\label{poteff}
\Ir_V (x) := V(x) -2\int \log |x-\xi|\!\ d\mu_V(\xi)
\end{align}
achieves its global minimum value on $(b^-, b^+) \setminus (\alpha^-, \alpha^+)$  only on the boundary of this set.
\end{itemize}
In the Hermite case, we have $V(x)= \frac{1}{2}x^2$ and the equilibrium measure $\mu_V$ is the semicircle law. In the Laguerre case, we may set $V(x)= \tau^{-1} x - (\tau^{-1}-1)\log(x)$ for $\tau\in(0,1]$ and $V(x)=+\infty$ for negative $x$. In this case, $\mu_V $ is the Marchenko-Pastur law $\operatorname{MP}(\tau)$. In both the Hermite and the Laguerre case, the assumption (A1), (A2) and (A3) are satisfied. We need one more definition related to outlying eigenvalues:
\begin{align}
\label{rate0}
\mathcal{F}_V^+(x) & = \begin{cases}
\mathcal{J}_V(x) - \inf_{\xi \in \R} \mathcal{J}_V(\xi) & \text{ if } \alpha^+\leq x \leq b^+, \\
\infty & \text{ otherwise, } 
\end{cases} \\ \label{rate0b}
\mathcal{F}_V^-(x) & = \begin{cases}
\mathcal{J_V}(x) - \inf_{\xi \in \R} \mathcal{J}_V(\xi) & \text{ if } b^-\leq x \leq \alpha^-, \\
\infty & \text{ otherwise. } 
\end{cases}
\end{align}
One may check that in the Hermite case, $\mathcal{F}_V^\pm=\mathcal{F}_H^\pm$ and in the Laguerre case, $\mathcal{F}_V^\pm=\mathcal{F}_L^\pm$, where $\mathcal{F}_H^\pm$ and $\mathcal{F}_L^\pm$ have been defined in the previous section.


\subsection{Basics on Large Deviations}
\label{sudular}
In order to be self-contained, let us recall the definition of large deviations  an some important tools (\cite{demboz98}).

\begin{defi}\label{defldp}
Let $E$ be a topological Hausdorff space and let $\mathcal{I}: E \rightarrow [0,\infty]$ be a lower semicontinuous function. 
We say that a sequence $(P_{n})_n$ of probability measures on $(E,\mathcal{B}(E))$ satisfies a large deviation principle (LDP) with 
rate function $\mathcal{I}$ and speed $a_n$ i f:
\begin{itemize}
\item[(i)] For all closed sets $F \subset E$:
\begin{align*}
\limsup_{n\rightarrow\infty} \frac{1}{a_n} \log P_{n}(F)\leq -\inf_{x\in F}\mathcal{I}(x)
\end{align*}
\item[(ii)] For all open sets $O \subset E$:
\begin{align*}
\liminf_{n\rightarrow\infty} \frac{1}{a_n} \log P_{n}(O)\geq -\inf_{x\in O}\mathcal{I}(x)\,.
\end{align*}
\end{itemize}
The rate function $\mathcal{I}$ is good if its level sets
$\{x\in E |\ \mathcal{I}(x)\leq a\}$ are compact for all $a\geq 0$.
We say that a sequence of $E$-valued random variables satisfies a LDP if their distributions satisfy a LDP. 
\end{defi}

\medskip

We will frequently use the following principle (\cite{demboz98},  p. 126).

\medskip

{\bf Contraction principle}. {\it Suppose that $(P_{n})_n$ satisfies an
LDP on $(E,\mathcal{B}(E))$ with good rate function $\mathcal{I}$ and speed $a_n$. Let
$f$ be a continuous mapping from $E$ to another topological Hausdorff space $F$.
Then $P_{n}\circ f^{-1}$ satisfies a LDP on $(F,\mathcal{B}(F))$
with speed $a_n$ and good rate function
\begin{align*}
\mathcal{I}'(y)=\inf_{\{ x \in E |f(x)=y\} } \mathcal{I}(x),\quad (y\in F).
\end{align*}}

\medskip

To prove our main large deviation principle, we will use 
a special extension of Baldi's theorem, which extends also Bryc's lemma. This new theorem is given in the Appendix.

 To apply this theorem  in our setting, we remark that the topological dual of $\mathcal{M}_p(T)$ is the space $\mathcal{C}_p(T)$ of bounded 
continuous  functions $f:T \to \mathcal H_p$ with the pairing
\begin{align*}
\langle \Sigma, f \rangle = \tr \int fd\Sigma .
\end{align*}

\medskip


\subsection{Large Deviations}
\label{suLDP}
\subsubsection{Random measures}

Our first LDP holds for $p\times p$ matrix measures 
\begin{align*}
\Sigma^{(N)} = \sum_{k=1}^N W_k \delta_{\lambda_k} ,
\end{align*}
whose support $(\lambda_1, \cdots, \lambda_N)$  
 is $\PNV$ distributed and where the distribution of weights $(W_1, \cdots, W_N)$ is $\mathbb{D}_{p,N}$ as in the case of classical ensembles. As explained in the introduction, this is precisely the distribution of the spectral measure of an $N\times N$ matrix $X_N$, drawn from the distribution 
 \eqref{matrixensemble}. Recall that under assumption (A1), the empirical measure of the eigenvalues converges to an equilibrium measure $\mu_V$, supported by $[\alpha^-,\alpha^+]$. The rate function of our large deviation principle involves the reference matrix measure  
\begin{align*}
\Sigma_V= \mu_V \cdot {\bf 1} \,.
\end{align*}
We recall that $\mathcal{F}_V^\pm$  has been defined in \eqref{rate0} and \eqref{rate0b}. The following theorem is the matrix counterpart of Theorem 3.1 in \cite{magicrules}. Note that in the scalar case, we had an additional parameter $\beta>0$, corresponding to the inverse temperature of the log-gas. In the matrix case, we choose to fix $\beta=2$ (for complex matrices) due to the nature of the matrix spaces.

\medskip

\begin{thm} \label{LDPgeneral}
Assume that the potential $V$ satisfies the assumptions (A1), (A2) and (A3). Then the sequence of spectral measures $\Sigma^{(N)}$ under $\PNV\otimes \mathbb{D}_{p,N}$ satisfies the LDP with speed $N$ and good rate function
\begin{align*}
\mathcal{I}_V(\Sigma) = \mathcal{K}(\Sigma_V\!\ |\!\ \Sigma) + 
\sum_{k=1}^{N^+} {\mathcal F}_{V}^+ (\lambda_k^+)  +  \sum_{k=1}^{N^-} {\mathcal F}_{V}^- (\lambda_k^-)
\end{align*}
if $\Sigma \in \mathcal{S}_{p,1}(\alpha^-,\alpha^+)$ and $\mathcal{I}_V(\Sigma) = \infty$ otherwise. 
\end{thm}

\medskip

\begin{rem} \label{LDPremark}
A natural extension of Theorem \ref{LDPgeneral} holds for potentials $V=V_N$ depending on $N$, provided that $V_N$ converges to a deterministic potential $V$ in an appropriate sense. For example, it holds if we suppose that $V_N:\mathbb{R} \rightarrow (-\infty,+\infty]$ is a sequence 
converging to $V$ uniformly on the level sets $\{V \leq M\}$, where $V$ satisfies assumptions (A1), (A2) and (A3)
and 
such that $V_N(x) \geq V(x)$.
\end{rem}

\subsubsection{Jacobi coefficients}
\label{subJac}

In the 
 cases of Hermite and Laguerre ensembles, the particular form of the distribution of the parameters $(A_1, B_1, \dots)$ and $(\zeta_1, \zeta_2, \cdots)$ respectively, allows us to 
 prove further LDP's for the spectral measure, independently of Theorem \ref{LDPgeneral}. They are in the subset $\mathcal{M}_{p,1,c}$ of compactly supported normalized matrix measures. Since we need a specific block structure, we assume $N=np$. 

\begin{thm} \label{LDPhermite}
Let $\Sigma^{(n)}$ be the spectral measure of $\frac{1}{\sqrt{np}} X_{n}$.  Assume that $X_{n}$ is distributed according to the Hermite ensemble $\operatorname{GUE}_N$ ($N=np$). Then the sequence $(\Sigma^{(n)})_n$ satisfies the LDP in $\mathcal{M}_{p,1,c}(\mathbb{R})$ with speed $pn$ and good rate function
\begin{align*}
\mathcal{I}_H(\Sigma) =  \sum_{k=1}^\infty  
\left[\tfrac{1}{2} \tr\!\  B_k^2  +  \tr\!\ G(A_kA_k^*)\right]
\end{align*}
where $B_k,A_k$ are the recursion coefficients of $\Sigma$ as in \eqref{polrek} if $\Sigma$ is non-trivial and $\mathcal I_H (\Sigma) = \infty$ if $\Sigma$ is trivial.
\end{thm} 

\medskip

\begin{thm} \label{LDPlaguerre}
Let $Y_n$ be distributed according to the Laguerre ensemble $\operatorname{LUE}_N(p\gamma_n)$, with 
($N=np$),
$\gamma_n\geq n$ an integer sequence such that $\frac{n}{\gamma_n}\to \tau \in (0,1]$ and let $\Sigma^{(n)}$ be the spectral measure of $\frac{1}{p\gamma_n}Y_n$ with a decomposition of recursion coefficients as in Section 2. Then the sequence $(\Sigma^{(n)})_n$ satisfies the LDP in $\mathcal{M}_{p,1,c}([0,\infty))$ with speed $pn$ and good rate function
\begin{align*}
\mathcal{I}_L(\Sigma) =   \sum_{k=1}^\infty  
\left[\tau^{-1}G(\zeta_{2k-1}) +  G(\tau^{-1}\zeta_{2k})\right],
\end{align*}
with $\zeta_k$ as in \eqref{zetadecomposition}. If $\Sigma$ is a trivial measure, we have $\mathcal{I}_L(\Sigma) = \infty$.
\end{thm}

\medskip

In order to prove LDP's for the spectral measures in terms of the recursion coefficients we need the following results for matrices of fixed size. The first and third are straightforward extensions of the scalar case, the second one can be found in \cite{GNRW}, with small changes to allow a general sequence of parameters.

\medskip

\begin{lem} \label{basicldps}
\begin{itemize}
\item[(i)] If $X \sim {\GUE}_p$ with $p$ fixed, then the sequence $(\tfrac{1}{\sqrt{n}} X)_n$ satisfies the LDP
 with speed $n$ and good rate function 
\begin{align*}\mathcal{I}_1(X) = \tfrac{1}{2} \operatorname{tr} X^2 .\end{align*}
\item[(ii)] Let $Y_n \sim \LUE_p(\gamma_n)$ with a positive sequence $(\gamma_n)_n$ such that $\frac{\gamma_n}{n}\to \gamma>0$, then the sequence $(\tfrac{1}{ n} Y_n)_n$ satisfies the LDP
 with speed $n$ and good rate function 
\begin{align*}
\mathcal{I}_2(Y) =   \gamma\!\ \tr\!\  G(\gamma^{-1}Y)
\end{align*}
if $Y$ is Hermitian and nonnegative and $\mathcal{I}_2(Y) = \infty$ otherwise.
\item[(iii)] Let $Z \sim \LUE_p(1)$ with $p$ fixed, that is, $Z=vv^*$ when $v$ is a vector of independent complex standard normal random variables. Then the sequence $(\tfrac{1}{n} Z)_n$ satisfies the LDP
 with speed $n$ and good rate function 
\begin{align*}
\mathcal{I}_3(Z) =  \tr\, Z 
\end{align*}
if $Z$ is Hermitian and nonnegative and $\mathcal{I}_3(Z) =\infty $ otherwise.
\end{itemize}
\end{lem}

\medskip


\section{Proof of Theorems \ref{sumruleHermite} and \ref{sumruleLag}: From large deviations to sum rules}
\label{sec:LDP2SR}

\textbf{Proof of Theorem \ref{sumruleHermite}:} \\
Consider the matrix measure $\Sigma^{(n)}$ with $N=np$ support points with density
\begin{align*}
\frac{1}{Z_V^N} \prod_{1\leq  i < j\leq N} |\lambda_i - \lambda_j|^2 \prod_{i=1}^N e^{- N\lambda_i^2/2} 
\end{align*}
and weight distribution $\mathbb{D}_{p,N}$ independent of the support points. By Theorem \ref{LDPgeneral}, the sequence $(\Sigma^{(n)})_n$ satisfies the LDP with speed $np$ and rate function $\mathcal{I}_V$, where $\mu_V$ is the semicircle law and furthermore $\mathcal{F}^\pm_V= \mathcal{F}_H^\pm$. That is, the rate function is precisely the left hand side of the equation in Theorem \ref{sumruleHermite}. On the other hand, $\Sigma^{(n)}$ is also the spectral measure of the random matrix $\frac{1}{\sqrt{pn}}X_n$, where $X_n \sim \operatorname{GUE}_{pn}$. By Theorem \ref{LDPhermite}, the sequence $(\Sigma^{(n)})_n$ satisfies also the LDP in the space of compactly supported measures with speed $np$ and rate function $\mathcal{I}_{H}$, the right hand side of the equation in Theorem \ref{sumruleHermite}. Since a large deviation rate function is unique, we must have $\mathcal{I}_V(\Sigma) = \mathcal{I}_H(\Sigma)$ for any compactly supported $\Sigma \in \mathcal{M}_{p,1,c}$. If $\Sigma$ is not compactly supported, it suffices to remark that the recursion coefficients cannot satisfy $\sup_n \left(||A_n|| + ||B_n||\right) <\infty$, as otherwise $J$ would be a bounded operator. But $\tr\!\ B^2 + \tr\!\ G(A)$ diverges as $||A||\to \infty$ or $||B||\to \infty$ and so the right hand side in Theorem \ref{sumruleHermite} equals $+\infty$.
\hfill $ \Box $ \\

\medskip

\textbf{Proof of Theorem \ref{sumruleLag}:} \\
Fix $\tau \in (0,1]$ and let $V(x) = \tau^{-1}x -(\tau^{-1}-1)\log x$ for $x\geq 0$ and $V(x)=+\infty$ if $x<0$. From Theorem \ref{LDPgeneral} we get that under the distribution $\mathbb{P}_N^V\otimes \mathbb{D}_{p,N}$ the sequence $(\Sigma^{(n)})_n$ satisfies the LDP with speed $N=np$ and rate function $\mathcal{I}_V$. In this case, the equilibrium measure is the Marchenko-Pastur law $\operatorname{MP}(\tau)$ multiplied by ${\bf 1}$. Further, we have $\mathcal{F}^\pm_V = \mathcal{F}^\pm_L$. So that,  $\mathcal{I}_V$ is nothing more than the left hand side of the sum rule in Theorem \ref{sumruleLag}. We would like to combine this result with the LDP in Theorem \ref{LDPlaguerre}, but since this requires integer parameters, we need to modify the potential slightly. Define $\gamma_n = \lceil n \tau^{-1} \rceil$ and consider the eigenvalue distribution with density
\begin{align} \label{eigenvaluevariation}
\frac{1}{Z_V^N} \prod_{1\leq  i < j\leq N} |\lambda_i - \lambda_j|^2 \prod_{i=1}^N \lambda_i^{p\gamma_n -pn} e^{-p\gamma_n \lambda_i} \mathbbm{1}_{[0,\infty)}(\lambda_i)
\end{align}
This is the eigenvalue distribution of the matrix $\frac{1}{p\gamma_n}Y_n$, when $Y_n \sim \operatorname{LUE}_{pn}(p\gamma_n)$. By Theorem \ref{LDPlaguerre}, the spectral measure of this matrix satisfies the LDP with speed $pn$ and rate function $\mathcal{I}_L$ which is the right hand side of the sum rule in Theorem \ref{sumruleLag}. We may as well write the eigenvalue density \eqref{eigenvaluevariation} as
\begin{align} \label{eigenvaluevariation2}
\frac{1}{Z_V^N} \prod_{1\leq  i < j\leq N} |\lambda_i - \lambda_j|^2 \prod_{i=1}^N e^{-pn V_n(\lambda_i)} ,
\end{align}
where 
\begin{align}
V_n(x) = \frac{\lceil n \tau^{-1} \rceil}{n} x - \left( \frac{\lceil n \tau^{-1} \rceil}{n} - 1 \right) \log x
\end{align}
for nonnegative $x$. Then $V_n(x)\geq V(x)$ for all $x$ and on the sets $\{x|\, V(x)\leq M\}$, the potentials $V_n$ converge uniformly to $V$. Note that the point 0 is included in the level sets of $V$ only if $\tau=1$. Therefore, by Remark \ref{LDPremark}, the spectral measure with support point density \eqref{eigenvaluevariation2} satisfies the same LDP as under the density $\mathbb{P}_{pn}^V$ and then with rate function $\mathcal{I}_V$. This yields $\mathcal{I}_V(\Sigma) = \mathcal{I}_H(\Sigma)$ for any compactly supported measure $\Sigma$. 
The extension to measures with non-compact support follows as in the proof of Theorem \ref{sumruleHermite}. 
\hfill $ \Box $ \\

\medskip


\section{Proof of Theorem \ref{LDPgeneral}: Spectral LDP for a general potential} 
\label{Sspopo}
This section is devoted to the proof of Theorem \ref{LDPgeneral}. We will follow the track of the proof developed for the scalar case in \cite{magicrules} and will often refer to this paper for more details. The main idea is to apply the projective method and study a family of matrix measures restricted to the support $I=[\alpha^-,\alpha^+]$ of the equilibrium measure and a fixed number of extremal eigenvalues. For  $\Sigma \in \mathcal{S}_{p}$ with
\begin{align} \label{SigmaS}
\Sigma = \Sigma_{|I} +  \sum_{i=1}^{N^+} \Gamma_i^+ \delta_{\lambda_i^+} + \sum_{i=1}^{N^-} \Gamma_i^- \delta_{\lambda_i^-}
\end{align}
we define the $j$-th projector $\pi_j$ by
\begin{align*}
\pi_j(\Sigma) = \Sigma_{|I} +  \sum_{i=1}^{N^+\wedge j} \Gamma_i^+ \delta_{\lambda_i^+} + \sum_{i=1}^{N^-\wedge j} \Gamma_i^- \delta_{\lambda_i^-} ,
\end{align*}
that is, all but the $j$-th largest and smallest eigenvalues outside of $I=[\alpha^-,\alpha^+]$ are omitted. Note that $\pi_j$ is not continuous in the weak topology. For this reason we need to change our topology on $\mathcal{S}_{p}$ by identifying $\Sigma$ as in \eqref{SigmaS} with the vector
\begin{align}\label{SigmaS2}
\big( \Sigma_{|I}, (\lambda_i^+)_{i \geq 1}, (\lambda^-_i)_{i\geq 1},(\Gamma^+_i)_{i\geq 1},(\Gamma^-_i)_{i\geq 1}\big)
\end{align}
 with $\lambda_i^+=\ap$ and $\Gamma_i^+=0$ if $i>N^+$ and $\lambda^-_i=\am$ and $\Gamma_i^-=0$ if $i>N^-$. Then we say that $\Sigma^{(n)}$ converges to $\Sigma$ if:
\begin{align} \label{strangetop}
\begin{split} \Sigma^{(n)}_{|I} \xrightarrow[n \rightarrow \infty ]{} \Sigma_{|I} & \text{ weakly and 
for every}\  i \geq 1 \\
\big( \lambda^+_i (\Sigma^{(n)}),\lambda^-_i(\Sigma^{(n)}),\Gamma^+_i(\Sigma^{(n)}),\Gamma^-_i&(\Sigma^{(n)})\big) \xrightarrow[n \rightarrow \infty ]{}\big( \lambda^+_i (\Sigma),\lambda^-_i(\Sigma),\Gamma^+_i(\Sigma),\Gamma^-_i(\Sigma)\big) . 
\end{split}
\end{align}
Analogously to Lemma 4.5 in \cite{magicrules}, one can show that on the smaller set $\mathcal{S}_{p,1}$ of normalized measures, this topology is (strictly) stronger than the weak topology. 

\medskip

Let for $j$ fixed and $N > 2j$ 
\[ \lambda^+(j) = ( \lambda_1^+, \dots, \lambda_j^+)\ , \ \lambda^-(j) = (\lambda_1^-, \dots, \lambda^-_{j})\,.\]
Then the following joint LDP holds for the largest and/or smallest eigenvalues, where we write $\mathbb R^{\uparrow j}$ (resp. $\mathbb R^{\downarrow j}$) for the subset of $\mathbb R^j$ of all vectors with non-decreasing (resp. non-increasing) entries and, with a slight abuse of notation, 
we write $\alpha^\pm$ for the vector $(\alpha^\pm,\dots ,\alpha^\pm)\in \mathbb{R}^j$.

\medskip

\begin{thm}
\label{LDPjextreme}
Let $j$ be a fixed integer and the potential $V$ such that (A1), (A2) and the control condition (A3) are satisfied.
\begin{enumerate}
\item If $b^-<\alpha^-$ and $\alpha^+<b^+$, then
the law of $(\lapj, \lambda^-(j))$ under $\PNV$ satisfies the LDP in $\R^{2j}$ with speed $N$ and good rate function
\begin{align*}
\mathcal{I}_{\lambda^\pm}(x^+, x^-)
 := 
\begin{cases}   
\sum_{k=1}^j 
\Fr_{V}^+(x_k^+)  + \sum_{k=1}^j 
\Fr_{V}^-(x^-_k)\
& \mbox{ if } 
(x_1^+, \dots , x_j^+)\in \mathbb R^{\downarrow j} \ \hbox{and}\  (x_1^-, \dots , x_j^-)\in \mathbb R^{\uparrow j}\\
 \infty & \mbox{ otherwise.}
\end{cases}
\end{align*}
\item If $b^-=\alpha^-$, but $\alpha^+<b^+$, the law of $\lapj$ satisfies the LDP with speed $N$ and good rate function
\begin{align*}
\mathcal{I}_{\lambda^+}(x^+)= 
\mathcal{I}_{\lambda^\pm}(x^+, \alpha^-)=
\begin{cases}   
\sum_{k=1}^j 
\Fr_{V}^+(x_k^+)
& \mbox{ if } 
(x_1^+, \dots , x_j^+)\in \mathbb R^{\downarrow j}\\
  \infty & \mbox{ otherwise.}
\end{cases}
\end{align*}
\item If $b^-<\alpha^-$, but $\alpha^+=b^+$, the law of $\lambda^-(j)$ satisfies the LDP with speed $N$ and good rate function
\begin{align*}
\mathcal{I}_{\lambda^-}( x^-) = \mathcal{I}_{\lambda^\pm}(\alpha^+, x^-) =
\begin{cases} 
\sum_{k=1}^j 
\Fr_{V}^-(x^-_k)& \mbox{ if }   (x_1^-, \dots , x_j^-)\in \mathbb R^{\uparrow j}\\
  \infty& \mbox{ otherwise.}
\end{cases}
\end{align*}
\end{enumerate}
\end{thm}

\medskip

\subsection{LDP for the restricted measure and extremal eigenvalues}

Suppose now that the distribution of $\Sigma^{(N)}$ is as in Theorem \ref{LDPgeneral} and the assumptions (A1), (A2) and (A3) are satisfied. By Lemma \ref{weightdistr}, we may decouple the weights and consider the (non-normalized) measure 
\begin{align}
\label{defsigmatilde}
\tilde\Sigma^{(N)} = \frac{1}{N} \sum_{k=1}^N v_kv_k^* \delta_{\lambda_k} ,
\end{align}
where the entries of $v_1,\dots ,v_N\in \mathbb{C}^p$ are independent complex standard normal distributed random vectors. 
 The original distribution can then be recovered as the pushforward under 
\begin{align*}
\tilde\Sigma \mapsto \tilde\Sigma(\mathbb{R})^{-1/2}\cdot \tilde\Sigma \cdot \tilde\Sigma(\mathbb{R})^{-1/2} .
\end{align*}  
Let $I(j):=I\setminus \{\lap_1,\lambda_1^-,\dots ,\lap_j,\lambda_j^- \}$ denote the interval $I$ without the $j$ largest and smallest eigenvalues. Analogously,  let $I^+(j):=I\setminus \{\lap_1,\dots ,\lap_j \}$ and $I^-(j):=I\setminus \{\lambda_1^-,\dots ,\lambda_j^- \}$. Then we write $\tilde\Sigma_{I(j)}^{(N)}$ for the restriction of $\tilde\Sigma^{(N)}$ to $I(j)$. We use the analogous notation for the restrictions to $I^+(j), I^-(j)$ and $I$. The main result in this section is a joint LDP for the restricted measure and the collection of extremal eigenvalues. 

\medskip

\begin{thm} \label{jointLDP} 
Suppose that the law of eigenvalues and weights is given by $\PNV \otimes \mathbb{G}_{p,N}$ and consider $\tilde\Sigma^{(N)}$ as a random element in $\mathcal{S}_p$ with topology \eqref{strangetop}. 
\begin{enumerate}
\item If $\bm<\am<\ap<\bp$, then for any fixed $j\in\N$, the sequence $\big( \SitnIj, \lapj, \lamj \big)$
satisfies the joint LDP with speed $N$ and good rate function
\begin{align*}
\mathcal{I}(\Sigma,x^+,x^-) =  \mathcal{K}(\Sigma_V\!\ | \!\  \Sigma) + \tr\!\ \Sigma (\R) -p+ \mathcal{I}_{\lambda^\pm}(x^+, x^-)
\end{align*}
\item If $\bm=\am$, but $\ap<\bp$ (or $\bp=\ap$, but $\am>\bm$), then, with the same notation as in Theorem \ref{LDPjextreme},
$
\big( \tilde\Sigma_{I^+(j)}^{(N)}, \lapj \big) (\text{or } \big( \tilde\Sigma_{I^-(j)}^{(N)}, \lamj  \big) \text{ respectively,} )
$
satisfies the LDP with speed $N$ and good rate function 
\begin{align*}
\mathcal{I}^+(\Sigma,x^+) = \mathcal{I}(\Sigma,x^+,\am) \quad (\text{or } \ \mathcal{I}^-(\Sigma,x^-) = \mathcal{I}(\Sigma,\ap,x^-) \text{ respectively})\,.
\end{align*}
\end{enumerate}
\end{thm}

\medskip

\proof
We show here only the first part of the theorem, for the other cases just omit the largest or smallest eigenvalues. 
To begin with, for $M>\max\{|\ap|,|\am|\}$, let $\lapjM$ (resp.$\lamjM$) be the collection of truncated eigenvalues 
\begin{align*}
\lambda_{M,i}^+ = \min\{ \lambda_i^+ , M\}\ \ \ (\text{resp.}\  \lambda_{M,i}^- = \max\{ \lambda_i^- , -M\})\,,
\end{align*}
for $i=1,\dots,j$. 
To further simplify the notations, set $\lambda_M^\pm(j):=(\lambda_{M,
1}^+,\dots,\lambda_{M,j}^+,\lambda_{M,1}^-,\dots ,\lambda_{M,j}^-)$. 
 The sequence $(\SitnIj, \lambda_M^\pm(j))$ is exponentially tight, since, with the compact set 
\begin{align*}
K_{\gamma,M}= \left\{ (\Sigma,\lambda)\in \mathcal{M}_p(I) \times \R^{2j} | \ \Sigma(I)\leq \gamma \cdot {\bf 1}, \lambda \in [-M,M]^{2j} \right\} ,
\end{align*}
we have 
\begin{align*}
\limsup_{N\to \infty} \frac{1}{N} \log  P\left( (\SitnIj, \lambda_M^\pm(j)) \notin K_{\gamma,M} \right) 
\leq \limsup_{N\to \infty} \frac{1}{N} \log P\left( \frac{1}{N} \sum_{k=1}^N v_kv_k^* > \gamma {\bf 1}\right) \leq - c\gamma ,
\end{align*}
where we used the fact that $\sum_{k=1}^N v_kv_k^*$ follows the $\operatorname{LUE}_p(N)$ distribution. We prove the joint LDP by applying Theorem \ref{newgeneral}. For this, let $D$ be the set of continuous $f: [\alpha^-,\alpha^+] \to \mathcal H_p$ such that for all $x\in [\alpha^-,\alpha^+]$, $f(x)<{\bf 1}$, i.e., the eigenvalues of $f(x)$ are smaller than 1. 
For $f\in D$ and $\varphi$ a bounded continuous function from $\mathbb R^{2j}$ to $\mathbb R$
, we consider the joint moment generating function 
\begin{align*}
{\mathcal{G}}_N(f, \varphi) = 
\E\left[ \exp\left\{ N \left(\tr \int f d\SitnIj +  \varphi(\lambda_M^\pm(j)) \right) \right\} \right] . 
\end{align*}
Since the weights $v_kv_k^*$ of $\SitnIj$ are independent, we may first integrate with respect to the $v_k$'s, so that 
\begin{align}
\label{condw}
{\mathcal{G}}_N(f,\varphi) &= \mathbb E \left[\exp\left(N \varphi(\lambda_M^\pm (j))\right) \prod_{i:\lambda_i \in I(j)} \mathbb E\left[\exp \left\{  \tr ( f(\lambda_i)v_kv_k^*)\right\} | \lambda_1, \dots, \lambda_n\right]\right]\\
&=   \mathbb E \left[\exp\left(N \varphi(\lambda_M^\pm (j))\right) \prod_{i:\lambda_i \in I(j)} \mathbb E\left[\exp \left\{  v_k^* f(\lambda_i)v_k\right\} | \lambda_1, \dots, \lambda_n\right]\right]\,.
\end{align}
Now, it is clear that
 for $v$ a standard normal complex vector in $\mathbb C^p$ and $A\in \mathcal H_p$ such that $A< 1$, we have
\begin{equation}
\label{defL}
\log  \mathbb{E} \left[ \exp \left( v^* A v \right) \right]  = - \log \det ({\bf 1}-A)= :L(A)\,
\end{equation}
so that \eqref{condw} becomes,
\begin{align*}
{\mathcal{G}}_N(f,\varphi) 
&= 
   \E\left[ \exp\left\{ N \left(\muunIj(L \circ f) + \varphi(\lambda_M^\pm(j))  \right)\right\} \right] \,,
\end{align*}
where $\muunIj$ is the restriction of the (scalar) empirical eigenvalue distribution to $I(j)$. It remains to calculate the expectation with respect to $\PNV$. 

By Theorem \ref{LDPjextreme}, the extremal eigenvalues $\lambda^\pm(j)$ of the spectral measure satisfy the LDP with speed $N$ rate function $\mathcal{I}_{\lambda^\pm}$. By the contraction principle (see \cite{demboz98} p.126), the truncated eigenvalues satisfy the LDP with rate function
\begin{align*}
\mathcal{I}_{M,\lambda^\pm}(x^\pm) = \begin{cases} \mathcal{I}_{\lambda^\pm}(x^+,x^-) & \text{ if } x^\pm = (x^+,x^-) \in [-M,M]^{2j}, \\ \infty & \text{ otherwise.} \end{cases}
\end{align*}
Since the truncated eigenvalues are bounded, we can conclude from 
 Varadhan's Integral Lemma (\cite{demboz98} p. 137) 
\begin{align} \label{jointLDPproof1}
\lim_{N\to \infty} \frac{1}{N}\log \E\left[ \exp\left\{ N \varphi(\lambda_M^\pm(j))  \right\} \right]
 = J(\varphi) := \sup_{y\in \mathbb{R}^{2j}} \{ \varphi (y) - \mathcal{I}_{M,\lambda^\pm}(y) \}\,.
\end{align}
Since $\muun$ satisfies a LDP with speed $N^2$, but we consider only the slower scale at speed $N$, we may replace it by the limit measure $\mu_V$ at a negligible cost. For the exact estimates, we may follow along the lines of \cite{magicrules} to conclude
\begin{align*}
\lim_{N\to \infty} \frac{1}{N} \log {\mathcal{G}}_{N}(f, \varphi) & = \lim_{N\to \infty} \frac{1}{N} \log \E\left[ \exp\left\{ N \left(\mu_V(L \circ f) + \varphi(\lambda_M^\pm(j))  \right)\right\} \right]\\
& = G(f) 
 +  J(\varphi) 
\end{align*}
where 
\[G(f) = \int L \circ f\, d\mu_V\]
and $L$ is given in \eqref{defL}. 
Theorem \ref{newgeneral} yields the LDP for $(\tilde\Sigma^{(n)}_I, \lambda^\pm_M)$ with good rate
\[\mathcal I(\Sigma, \lambda) = G^*(\Sigma) + \mathcal I_{M, \lambda^\pm}(\lambda)\,,\]
once we show the second assumption therein is satisfied. 

Theorem 5 of \cite{rockafellar2} identifies $G^*$ as
\begin{align}
\label{rock}
G^*(\Sigma) =  \int L^*(h) d\mu_V + \int r\left(\frac{d\Sigma^S}{d\theta}\right) d\theta ,
\end{align}
where:
\begin{itemize}
\item $L^*$ is the
convex dual of $L$ 
\item $r$ its recession function
\item   the Lebesgue-decomposition of $\Sigma$ with respect to $\mu_V$ is
\[d\Sigma(x) = h(x)d\mu_V(x) + d\Sigma^S(x)\] 
\item
$\theta$ is any scalar measure such that $\Sigma^S$ is absolutely continuous with respect to $\theta$.
\end{itemize} 
We begin by calculating $L^*$ and $r$. By definition,
\begin{align*}
L^*(X) =& \sup_{Y\in \mathcal H_p} \left\{ \tr(XY)- L(Y) \right\} .
\end{align*} 
 The recession function is
\begin{align*}
r(X) = \sup \left\{ \tr(XW) \left| \ ||W||_\infty < 1 \right. \right\} .
\end{align*}

The function $L$ is convex (as in the scalar case, apply H\"{o}lder's inequality in the definition (\ref{defL})) and analytic. The supremum is then reached at a critical value. We denote by $\mathcal{D}[F(Y)]$ the Fr\'{e}chet derivative of a function $F:\mathcal H_p\to \mathbb{R}$ at $Y$ and look for $Y$ such that
\begin{equation}\label{critical}\mathcal{D}[ \tr (XY) - L(Y)](Z) = 0\end{equation}
for every $Z$. It is well known that, as functions of $Y$ for $X$ fixed, $\mathcal{D}[\tr(XY)] (Z) = \tr (XZ)$ and $\mathcal{D}[\log \det Y](Z) = \tr (Y^{-1}Z)$ so that (\ref{critical}) becomes, by the chain rule,
\[\tr (XZ) - \tr (({\bf 1} - Y)^{-1}Z) = 0\]
for every $Z$ i.e.
$X - (I - Y)^{-1} = 0$ hence $Y= {\bf 1} - X^{-1}$ and 
\begin{equation}
\label{Lstar}
L^*(X) =  \tr \left(X - {\bf 1}\right) + \log\det (X^{-1}) =  \tr X - p + \log\det (X^{-1})=\tr\!\ G(X)\,.\end{equation}

If $X$ has a negative eigenvalue, then $r(X) = \infty$. For $X$ nonnegative definite, the supremum is attained for $W = {\bf 1}$, so that 
\begin{equation}\label{rec}r(X) =  \tr X . \end{equation}
 Gathering (\ref{Lstar}) and (\ref{rec}) and plugging into  (\ref{rock}) we get
\begin{align*}
G^*(\Sigma) & =   \tr \int h d\mu_V - \int \log \det h\!\  d\mu_V  -p +  \tr \int  d\Sigma^S \\
& =  \mathcal{K}(\Sigma_V | \Sigma ) +   \tr\!\ \Sigma (I) -p .
\end{align*}
It remains to show that $G^*$ is sufficiently convex.  
A measure $\Sigma \in \mathcal{M}_p([\alpha^-,\alpha^+])$ is a point of strict convexity (an exposed point) of $G^*$ if there exists an exposing hyperplane $f \in C_p([\alpha^-,\alpha^+])$, such that
\begin{align}\label{exposedpoint}
\tr \int fd\Sigma - G^*(\Sigma) > \tr \int fd\zeta - G^*(\zeta)
\end{align} 
for all $\zeta \not = \Sigma$ (see \eqref{exposinghyper}). Let $\Sigma = h_\Sigma \mu_V$ be absolutely continuous with respect to $\mu_V$ with positive definite continuous density $h_\Sigma$ and choose
\begin{align*}
f=  {\bf 1} - h_\Sigma^{-1}\,.
\end{align*}
Then $f\in D$ and, by continuity and compactness of $[\alpha^-,\alpha^+]$, there exists a $\gamma >1$ such that $\gamma f \in D$. 
Let $\zeta = h_\zeta \mu_V + \zeta^S$ the Lebesgue-decomposition of $\zeta$.
Recalling the 
representation (\ref{rock}) and (\ref{Lstar}),
inequality \eqref{exposedpoint} is 
satisfied as soon as
\begin{align}
\label{CS1}\int \log \det h_\Sigma\!\ d\mu_V - p
 > \int \log\det h_\zeta\!\   d\mu_V- \tr \int h_\Sigma^{-1} d\zeta \,.\end{align}
Since $\tr \int h_\Sigma^{-1} d\zeta \geq \tr \int h_\Sigma^{-1}h_\zeta d\mu_V$, it is enough to prove
\begin{align}
\label{CS2}\int \log \det h_\Sigma\!\ d\mu_V - p
 > \int \log\det h_\zeta\!\   d\mu_V- \tr \int h_\Sigma^{-1} h_\zeta d\mu_V\,.
\end{align} 
This inequality follows from
\begin{align}\label{logineq}
\log \det A -  \log \det B > p - \tr (A^{-1}B)
\end{align}
for Hermitian positive $A\neq B$. 
In order to prove \eqref{logineq}, write
\begin{align}\label{logineq2}
\log \det A -  \log \det B = \sum_{i=1}^p \log \left( \lambda_i(A) \lambda_i(B)^{-1} \right) \geq 
p - \sum_{i=1}^p \lambda_i(A^{-1}) \lambda_i(B) 
\end{align}
with $\lambda_i(A), \lambda_i(B)$ the eigenvalues of $A,B$ written in any order. If we choose to order the eigenvalues of $A^{-1}$ in decreasing order (i.e. those of $A$ increasing) and those of $B$ in increasing order, it follows from the Hardy-Littlewood rearrangement inequality (see \cite{mirsky}) that
\begin{align*}
\sum_{i=1}^p \lambda_i(A^{-1}) \lambda_i(B) \leq \tr (A^{-1}B) .
\end{align*}
With this ordering of eigenvalues, \eqref{logineq2} is strict unless $A,B$ have the same eigenvalues. If all eigenvalues of $A$ and $B$ coincide, then the left hand side of \eqref{logineq} is 0, while the right hand side is $p-\tr H$ with $\det H=1$. The minimum value of $\tr H$ is $p$ which is achieved only for $H={\bf 1}$, in which case $A=B$. We get that $\Lambda^*$ is strictly convex at all points $\Sigma = h\mu_V$ with $h$ positive definite and continuous.

It remains to show that the set of exposed points of $G^*$ is dense in $\mathcal{M}_p([\alpha^-,\alpha^+])$. For a given $\Sigma \in \mathcal{M}_p([\alpha^-,\alpha^+])$, we divide $[\alpha^-,\alpha^+]$ by dyadic points into intervals \[I_{k,n} = [\alpha^-+ (k-1)(\alpha^+-\alpha^-)/2^n,\alpha^-+k(\alpha^+-\alpha^-)/2^n]\] and put $\overline{I}_{k,n}= [\alpha^-+ (k-1)(\alpha^+-\alpha^-)/2^n+ 2^{-2n},\alpha^-+k(\alpha^+-\alpha^-)/2^n-2^{-2n}]$, i.e., $\overline{I}_{k,n}$ is constructed from ${I}_{k,n}$ by cutting off subintervals of length $2^{-2n}$. Define $\overline{h}_n$ on $\overline{I}_{k,n}$ as 
\begin{align*}
{\overline{h}_n}_{\big | \overline{I}_{k,n}} \equiv \left( (1-2^{-2n})\Sigma(I_{k,n}) + 2^{-2n} \cdot I \right)\Sigma_V(I_{k,n})^{-1} 
\end{align*}
and let $h_n$ be the continuous function on $[\alpha^-,\alpha^+]$ obtained by linear interpolation of the step function $\overline{h}_n$. Then $h_n$ is positive definite and continuous on $[\alpha^-,\alpha^+]$ and as in the scalar case, $h_n\cdot \Sigma_V$ converges weakly to $\Sigma$. This concludes the proof of the LDP for $\big( \SitnIj, \lambda_M^\pm(j) \big)$.\\

In order to extend the LDP to the untruncated eigenvalues, note that the LDP for $(\lapj, \lamj )$ implies the exponential tightness of the (unrestricted) extremal eigenvalues  
that is, for every $K>0$ there exists a $M$ such that
\begin{align*}
\limsup_{N\to \infty} \frac{1}{N} \log P\big( \lambda_1^+>M \text{ or } \lambda_1^-<-M \big) \leq -K .
\end{align*}
In particular, 
\begin{align*}
\lim_{M\to \infty}  \limsup_{N\to \infty} \frac{1}{N} \log P\big( \lambda_M^\pm(j) \neq \lambda^\pm(j) \big) = -\infty ,
\end{align*}
so that as $M\to \infty$, the truncated eigenvalues are exponentially good approximation of the unrestricted ones. Moreover, $(\SitnIj,\lambda_M^\pm(j) )$ are exponentially good approximations of $(\SitnIj,\lambda^\pm(j) ) $. 
By Theorem 4.2.16 in \cite{demboz98} $(\SitnIj,\lambda^\pm(j) ) $ satisfies the LDP with speed $n$ and rate function 
\begin{align*}
\mathcal{I}(\Sigma,x^\pm) & = \mathcal{K}(\Sigma_V\!\ |\!\ \Sigma) +  \tr\!\ \Sigma (I) -p + \mathcal{I}_{\lambda^\pm}(x^\pm) \\ 
& = \mathcal{K}(\Sigma_V\!\ |\!\ \Sigma) +   \tr\!\ \Sigma (I) -p +  \sum_{i=1}^ {j} \Fr^+(x^+_i) + \Fr^-(x^-_i)\,,
\end{align*}
which ends the proof of Theorem \ref{jointLDP}. \hfill $ \Box $ \\


\subsection{LDP for the projective family}

\begin{thm}\label{projLDP}
For any fixed $j$, the sequence of projected spectral measures $\pi_j(\Sitn)$ as elements of $\mathcal{S}_p$ with topology \eqref{strangetop} satisfies the LDP with speed $N$ and rate function
\begin{align*}
\tilde{\mathcal{I}}_j(\tilde{\Sigma}) = \mathcal{K}(\Sigma_V \!\ |\!\ \tilde{\Sigma}) + \tr\!\ \tilde\Sigma (I) -p+ \sum_{i=1}^{N^+\wedge j} \left(\Fr^+_V(\lambda^+_i)+\tr\!\  \Gamma_i^+\right) + \sum_{i=1}^{N^-\wedge j}\left(\Fr^-_V(\lambda^-_i)+\tr\!\ \Gamma_i^-\right).
\end{align*}
\end{thm}

\medskip

\proof
The proof is similar to the proof of Theorem 4.3 in \cite{magicrules} and we omit the details for the sake of brevity. It is a combination of the LDP in Theorem \ref{jointLDP} and the LDP of the weights $\frac{1}{N} \Gamma_k=\frac{1}{N}v_kv_k^*$ corresponding to the extreme eigenvalues. Indeed $\Gamma_i \sim \LUE_p(1)$, so by Proposition \ref{basicldps} (iii), each individual weight $\frac{1}{N}\Gamma_k$ satisfies the LDP with speed $N$ and rate function 
\begin{align*}
\mathcal{I}_3 (X) = \begin{cases} 
 \tr\!\  X & \text{ if } X\geq 0, \\
\infty & \text{ otherwise.}
\end{cases}
\end{align*} 
The independence of the weights and an application of the contraction principle complete then the proof. 
\hfill $ \Box $ \\

\medskip

In order to come back to a normalized matrix measure in $\mathcal{S}_{p,1}$, we note that the LDP for $\pi_j(\Sitn)$ also implies the joint LDP for
\begin{align*}
\big(\pi_j(\Sitn) , \pi_j(\Sitn)(\R) \big) ,
\end{align*}
with the rate function
\begin{align*}
\overline{\mathcal{I}}_j(\tilde\Sigma, Z) = \tilde{\mathcal{I}}_j(\tilde{\Sigma})
\end{align*}
if $\tilde\Sigma(\mathbb{R})=Z$ and $\overline{\mathcal{I}}_j(\tilde\Sigma, Z)=\infty$ otherwise.
Keeping the weights along the way, we may apply the projective method (the Dawson-G\"artner Theorem, p. 162 in the book of \cite{demboz98}) to the family of projections $(\pi_j(\Sitn) , \pi_j(\Sitn)(\R) )_j$ and get a LDP for the pair 
$(\Sitn,\Sitn(\R))$ with rate function
\begin{align*}
\overline {\mathcal{I}}(\tilde\Sigma,Z) = \sup_j \overline{\mathcal{I}}_j(\tilde\Sigma, Z)\,.
\end{align*}
This rate function equals $+\infty$ unless $\tilde\Sigma(\mathbb{R})=Z$ and in this case is given by
\begin{align} \label{ratepair}
\overline {\mathcal{I}}(\tilde\Sigma,Z)
= \mathcal{K}(\Sigma_V\ | \ \tilde\Sigma) + \tr Z- p+\sum_{i=1}^{N^+} \Fr(\lambda^+_i)+\sum_{i=1}^{N^-}\Fr(\lambda^-_i) .
\end{align}
We remark that normalizing the matrix measure is not possible unless we keep track of the total mass for any $j$, as the mapping $\tilde\Sigma\mapsto \tilde\Sigma(\R)^{-1/2} \tilde\Sigma \tilde\Sigma(\R)^{-1/2}$ is not continuous in the topology \eqref{strangetop}. 
However, we may now apply the continuous mapping $(\tilde\Sigma,Z) \mapsto Z^{-1/2} \tilde\Sigma Z^{-1/2}$ and obtain a LDP for the sequence of measures $\Sigma^{(N)}$ in $\mathcal{S}_{p,1}$. The rate function is 
\begin{align*}
\mathcal{I}(\Sigma) = \inf_{\tilde\Sigma= Z^{1/2} \Sigma Z^{1/2},\, Z>0} \tilde{\mathcal{I}}(\tilde\Sigma) = \inf_{Z>0} \tilde{\mathcal{I}}(Z^{1/2} \Sigma Z^{1/2}) .
\end{align*}
By \eqref{ratepair}, we need to minimize over positive definite $Z \in \mathcal H_p$ the function
\begin{align*}
& \quad   - \int\log\det \left(\frac{d(Z^{1/2} \Sigma Z^{1/2})}{d\mu_V}\right)d\mu_V +  \tr Z -p \\
& =  - \int\log\det \left(Z^{1/2}\frac{d \Sigma }{d\mu_V}Z^{1/2}\right)d\mu_V +  \tr Z -p\\
& = - \int\log\det \left(\frac{d \Sigma }{d\mu_V}\right)d\mu_V - \log\det Z +  \tr Z -p
\\ &=  - \int\log\det \left(\frac{d \Sigma }{d\mu_V}\right)d\mu_V + \mathcal I_2 (Z)\,.
\end{align*}
The term 
 $\mathcal{I}_2(Z)$ comes from Lemma \ref{LDPgeneral} (ii) with $\gamma =1$ and attains its minimal value $0$ for $Z = {\bf 1}$.  \hfill $ \Box $ \\

We have obtained the LDP claimed in Theorem \ref{LDPgeneral} on the subset $\mathcal{S}_{p,1}$ in the topology induced by \eqref{strangetop}. On $\mathcal{S}_{p,1}$ this is stronger than the weak topology  and the rate function can be extended to $\mathcal{M}_{p,1}$ by setting $\mathcal{I}(\Sigma)=\infty$ for $\Sigma \notin \mathcal{S}_{p,1}$. This yields Theorem \ref{LDPgeneral}.

\subsection{Proof of Remark \ref{LDPremark}}

Let $A$ be a measurable subset of $\mathcal{M}_{p,1}$ and set 
\begin{align*}
A_N = \left\{ (\lambda,W) \in \mathbb{R}^N \times \mathcal H_p^N \left|\, \sum_{k=1}^N W_k \delta_{\lambda_k} \in A \right. \right\} .
\end{align*}
The LDP for $\Sigma^{(N)}$ with eigenvalue distribution $\mathbb{P}_N^{V_N}$ will follow from the LDP for eigenvalue distribution $\mathbb{P}_N^{V}$ once we show 
\begin{align} \label{LDPequiv1}
\limsup_{N\to \infty} \frac{1}{N} \log (\mathbb{P}_N^{V_N}\otimes \mathbb{D}_{p,N})(A_N) \leq \limsup_{N\to \infty} \frac{1}{N} \log (\mathbb{P}_N^{V}\otimes \mathbb{D}_{p,N})(A_N)
\end{align}
and
\begin{align} \label{LDPequiv2}
\liminf_{N\to \infty} \frac{1}{N} \log (\mathbb{P}_N^{V_N}\otimes \mathbb{D}_{p,N})(A_N) \geq \liminf_{N\to \infty} \frac{1}{N} \log (\mathbb{P}_N^{V}\otimes \mathbb{D}_{p,N})(A_N) .
\end{align}
In fact, this does not require $A$ to be closed or open, respectively. 
For this, let
\begin{align*}
\Gamma_N^V(A_N) = \iint_{A_N} \prod_{1\leq  i < j\leq N} |\lambda_i - \lambda_j|^2  \prod_{i=1}^{N} e^{-NV(\lambda_i)} d\lambda \!\ d\mathbb{D}_{p,N}(W)\,, 
\end{align*}
and define $\Gamma_N^{V_N}(A_N)$ analogously, with $V$ replaced by $V_N$. Since $V_N \geq V$, we have
\begin{align}\label{gammaineq1}
\Gamma_N^{V_N}(A_N) \leq \Gamma_N^V(A_N).
\end{align}
To get a reverse inequality, let $K_{N,M}$ be the set of $(\lambda,W) \in \mathbb{R}^N \times \mathcal H_p^N$, where $V(\lambda_i) \leq M$ for all $i$. Then
\begin{align*}
\Gamma_N^{V_N}(A_N) \geq \Gamma_N^{V_N}(A_N\cap K_{N,M}) \geq \left( \inf_{x : V(x) \leq M } e^{V(x)-V_N(x)} \right)^N \Gamma_N^{V}(A_N\cap K_{N,M}) .
\end{align*}
Since by assumption $e^{V(x)-V_N(x)}$ converges to 1 uniformly on $\{x|\, V(x) \leq M\}$, this implies
\begin{align} \label{gammaineq2} 
\lim_{N\to \infty} \frac{1}{N} \log \frac{\Gamma_N^{V_N}(A_N)}{\Gamma_N^{V}(A_N)} \geq \lim_{N\to \infty} \frac{1}{N} \log \frac{\Gamma_N^{V}(A_N\cap K_{N,M})}{\Gamma_N^{V}(A_N)} .
\end{align}
If we take now $A=\mathcal{M}_{p,1}$, then $\Gamma_N^{V}(A_N)=Z^V_N$ and the right hand side of \eqref{gammaineq2} becomes
\begin{align*}
\lim_{N\to \infty} \frac{1}{N} \log \mathbb{P}_N^V ( \forall\, i: \, V(\lambda_i) \leq M ) . 
\end{align*}
By the LDP for the extreme eigenvalues, Theorem \ref{LDPjextreme}, this limit tends to 0 as $M\to \infty$.   Together with \eqref{gammaineq1}, we have shown that for $A=\mathcal{M}_{p,1}$
\begin{align*}
\lim_{N\to \infty} \frac{1}{N} \log \frac{\Gamma_N^{V_N}(A_N)}{\Gamma_N^{V}(A_N)} = \lim_{N\to \infty} \frac{1}{N} \log \frac{Z_N^{V_N}}{Z_N^V} = 0 .
\end{align*}
Since $(\mathbb{P}_N^{V_N}\otimes \mathbb{D}_{p,N})(A_N) = (Z_N^V)^{-1}\Gamma_N^V(A_N)$, the inequality \eqref{gammaineq1} leads to the inequality \eqref{LDPequiv1} and the arguments for \eqref{gammaineq2} yield
\begin{align*}
\liminf_{n\to \infty} \frac{1}{N} \log (\mathbb{P}_N^{V_N}\otimes \mathbb{D}_{p,N})(A_N) \geq \liminf_{N\to \infty} \frac{1}{N} \log (\mathbb{P}_N^{V}\otimes \mathbb{D}_{p,N})(A_N\cap K_{N,M})
\end{align*}
for any $M\geq 0$. Letting $M\to \infty$, this implies inequality \eqref{LDPequiv2}, as by the LDP for the extreme eigenvalues we have 
\begin{align*}
\lim_{M\to \infty} \limsup_{N\to \infty} \frac{1}{N} \log (\mathbb{P}_N^{V}\otimes \mathbb{D}_{p,N})( K_{N,M}^c) = -\infty .
\end{align*}
 \hfill $ \Box $ \\


\section{Proof of Theorems \ref{LDPhermite} and \ref{LDPlaguerre}}
\label{sec6}
\subsection{Hermite block case}
\label{sHerm}
The starting point for the proof of Theorem \ref{LDPhermite} is the following block-tridiagonal representation of the Gaussian ensemble. It is a matrix extension of a famous result of Dumitriu and Edelman \cite{dumede2002}.

\medskip

\begin{lem} \label{triG}
Let $D_k \sim \GUE_p$ and $C_k$  be Hermitian non-negative definite
 such that $C_k^2 \sim  \LUE_p(p(n-k))$ for $k=1,\dots ,n$ and let all these matrices be independent. Then the $p\times p$ 
 spectral measure of the block-tridiagonal matrix 
\begin{align} \label{triGmatrix} \mathcal{G}_n = 
 \begin{pmatrix} 
  D_1 & C_1    &         &         \\
                C_1 & D_2    & \ddots  &         \\
                    & \ddots & \ddots  & C_{n-1} \\
                    &        & C_{n-1} & D_n
\end{pmatrix}
\end{align}
has the same distribution as the spectral measure of the Hermite ensemble $\GUE_{pn}$.
\end{lem}

\medskip

\proof Starting from a matrix $X_n$ distributed according to the Hermite ensemble $\GUE_{pn}$, we can construct the tridiagonal matrix $\mathcal{G}_n$ as
\begin{align*}
\mathcal{G}_n = TX_nT^*,
\end{align*}
where $T$ is unitary and leaves invariant the subspaces spanned by the first $p$ unit vectors. Consequently, the spectral measure of $X_n$ is also the spectral measure of $\mathcal{G}_n$. The transformation $T$ is a composition of unitaries $T_1,\dots ,T_{n-1}$ analogous to the ones used by \cite{dumede2002}. To illustrate the procedure, we construct the first transformation $T_1$. 
By $x_{i,j}$ we denotes the $p\times p$ subblock of $X_n$ in position $i,j$, let $\bar{x}_1= (x_{1,2},\dots ,x_{1,n})^*$ and $\bar{X}= (x_{i,j})_{2\leq i,j\leq n}$. With this notation, $X_n$ can be structured as
\begin{align*}
X_n = \begin{pmatrix}
x_{1,1} & \bar{x}_1^* \\
\bar{x}_1 & \bar{X}
\end{pmatrix} .
\end{align*}
Note that the Gaussian distribution implies that all (square) blocks are almost surely invertible. Then, set
\begin{align*}
\xi & = 
[(x_{2,1}^*)^{-1} (\bar{x}_1^*\bar{x}_1)(x_{2,1})^{-1}]^{1/2} x_{2,1}
\in \mathbb{M}_{p,p} \\
\Gamma & = (\xi^*,  {\bf 0}, \dots,  {\bf 0})^* \in \mathbb{M}_{(n-1)p,p}
\end{align*}
and define for $Z\in \mathbb{M}_{(n-1)p,p}$ the block-Householder reflection
\begin{align*}
H(Z) = I_{(n-1)p} - 
2 Z(Z^*Z)^{-1} Z^* .
\end{align*}
If we set $Z=\Gamma-\bar{x}_1$ one may check that
\begin{align*}
\Gamma^*\Gamma = \xi^*\xi = \bar{x}_1^* \bar{x}_1, \quad \Gamma^*\bar{x}_1 = \xi^* x_{2,1} = x_{2,1}^* \xi= \bar{x}_1^*\Gamma, \quad 
Z^*Z = -2Z^*\bar x_1
\end{align*}
and  
$H(Z) \bar x_1 = \Gamma$.
We extend $H(Z)$ to an operator $\hat{H}$ on $\mathbb{C}^{np}$ leaving the first $p$ unit vectors invariant, which yields
\begin{align*}
\hat{H} X_n \hat{H}^* = \begin{pmatrix} {\bf 1} & 0 \\
0 & H(Z) \end{pmatrix}
\begin{pmatrix} x_{1,1} & \bar{x}_1^* \\
\bar{x}_1 & \bar{X} \end{pmatrix}
\begin{pmatrix} {\bf 1} & 0\\
0 & H(Z)^* \end{pmatrix}
= \begin{pmatrix} x_{1,1} & \Gamma^* \\
\Gamma & H(Z)\bar{X} H(Z)^* \end{pmatrix} .
\end{align*}
Finally, let $W\in \mathbb{M}_{(n-1)p,(n-1)p}$ be the unitary block-diagonal matrix with the blocks $((\xi^*\xi)^{1/2}\xi^{-1},{\bf 1},\dots ,{\bf 1})$ on the diagonal and extend $W$ to an operator $\hat{W}$ on $\mathbb{C}^{np}$ as we did with $H(Z)$. Then $T_1=\hat{W}\hat{H}$ is unitary, leaves the subspace spanned by the first $p$ unit vectors invariant and
\begin{align*}
T_1 X_n T_1^* = 
\begin{pmatrix} x_{1,1} & \tilde{\Gamma}  \\
\tilde{\Gamma}^* & WH(Z)\bar{X}H(Z)^*W^* \end{pmatrix} .
\end{align*}
with 
\begin{align*}
\tilde{\Gamma}= ((\bar{x}_1^*\bar{x}_1)^{1/2}, {\bf 0},\dots ,  {\bf 0}) .
\end{align*}
The block $\bar{x}_1^*\bar{x}_1$ is distributed according to 
$\LUE_p (p(n-1))$ and since the definition of $W$ and $H(Z)$ is independent of $\bar{X}$, the block $WH(Z)\bar{X}H(Z)^*W^*$ is again a matrix of the Gaussian ensemble $\GUE_{p(n-1)}$. The assertion follows then from an iteration of these reflections.\hfill $ \Box $ \\
\medskip

\textbf{Proof of Theorem \ref{LDPhermite}:} \\
By Lemma \ref{triG}, the spectral measure $\Sigma^{(n)}$ is also the spectral measure of the rescaled matrix $\frac{1}{\sqrt{np}} \mathcal{G}_n$. If we consider each block entry of this matrix separately, we are up to a linear change of the speed in the setting of Lemma \ref{basicldps}. Thus, for any fixed $k$, the block $D_k^{(n)} := D_k/\sqrt{np}$ of the matrix  in \eqref{triGmatrix} satisfies the LDP in $\mathcal H_p$ with speed $pn$ and rate function $\mathcal I_1$. 
Similarly, if we let $C_k^{(n)}= C_k/{\sqrt{np}}$, then the 
 block $(C_k^{(n)})^2$ satisfies the LDP with speed $pn$ and rate function $\mathcal I_2$ or equivalently, $C_k^{(n)}$ 
satisfies the LDP with speed $np$ and good rate 
\begin{align*}
 \mathcal{I}'_2(Y)  
 = \mathcal{I}_2(Y^2) 
\end{align*}
if $Y$ is nonnegative definite and $ \mathcal{I}'_2(Y)=\infty$ otherwise. Since the block entries are independent, we get a joint LDP for any fixed collection $(D_1^{(n)},C_1^{(n)} ,\dots ,D_k^{(n)})$ with rate 
given by the corresponding sum of rate functions $\mathcal{I}_1$ and $\mathcal{I}'_2$.

Now, we follow the strategy developed in \cite{gamboa2011large} for the scalar case.
The random matrix measure $\Sigma^{(n)}$ belongs to $\mathcal M_{p, c}^1(\mathbb{R})$. Since the topology $\mathcal T_m$ of the convergence of moments on 
$\mathcal M_{p, c}^1(\mathbb{R})$ is stronger than the trace $\mathcal T_w$ of the weak topology, it is enough to prove the LDP with respect to $\mathcal T_w$.

For each $k > 0$, the subset $X_k$ of matrix probability measures with support in $[-k, k]$ is compact for $\mathcal T_m$.
Since the extremal eigenvalues satisfy the LDP with speed $N$ and a rate function tending to infinity, we deduce that $\Sigma^{(N)}$ is exponentially tight in $\mathcal T_m$.

The mapping 
\begin{align*}
m: \mathcal M_{p, c}^1(\mathbb{R}) \to \mathcal H_p^{\mathbb N_0}, \qquad m(\Sigma) : = \left(m_k(\Sigma) := \int x^k d\Sigma(x)\right)_{k \geq 1} 
\end{align*}
being a continuous injection, the LDP of $\Sigma^{(N)}$ in $\mathcal M_{p, c}^1(\mathbb{R})$ is then a consequence of the following LDP on the sequence of moments and of the inverse contraction principle (see \cite{demboz98} Theorem 4.2.4 and the subsequent Remark (a)). 

\begin{prop} 
The sequence $(m( \Sigma^{(n)})_n)$ satisfies the LDP in $\mathcal H_p^{\mathbb N_0}$ with speed $np$ and good rate function $\mathcal I_{\operatorname{m}}$ defined as follows. This function is finite in $(m_1, m_2, \dots)$ if and only if there exists a sequence $(B_1, A_1, \dots)\in \mathcal H_p^{\mathbb N_0}$ with $A_k > 0$, such that
\[\sum_{k=1}^\infty \frac{1}{2}\tr\!\  B_k^2  + \tr G(A_k^2) < \infty\]
and such that
\begin{equation} \label{momrelation}
(m_r)_{i,j} = \langle e_i J^r e_j\rangle \ , \ i,j = 1, \dots, p , \ r \geq 1
\end{equation} 
where $J$ is the infinite block Jacobi matrix with blocks $(B_1, A_1, \dots)$ as in \eqref{jacobimatrix}.

In that case
\begin{equation}
\mathcal I_{\operatorname{m}} = \sum_{k=1}^\infty \frac{1}{2}\tr B_k^2  + \tr G(A_k^2)\,.
\end{equation}
\end{prop}
\proof First, as we said in the beginning of this proof, for fixed $k$, $(D_1^{(n)}, C_1^{(n)}, \dots, D_k^{(n)} )$ satisfies the LDP in $\mathcal{H}_p^{2k-1}$ with speed $np$ and good rate function
\[\mathcal{I}^{(k)}(D_1,C_1,\dots , D_{k}) =  \sum_{j=1}^k \frac{1}{2}\tr D_j^2  + \sum_{j=1}^{k-1}\tr G(C_k^2)\,.\]
If $J$ is the $kp\times kp$ Jacobi matrix build from the blocks $D_1,C_1,\dots ,D_k$, then the moments $(m_1(\Sigma^{(n)}), \dots, m_{2k -1}(\Sigma^{(n)}))$ of the spectral measure of $J$ are given by \eqref{momrelation} and depend continuously on $D_j,C_j$. By the contraction principle, the sequence $(m_1(\Sigma^{(n)}), \dots, m_{2k -1}(\Sigma^{(n)}))$ satisfies the LDP with speed $np$ and good rate function ${\mathcal{I}}_{\operatorname{m}}^{(k)}$ defined as follows. It is infinite in $(m_1,\dots ,m_{2k-1})$ unless there exist block coefficients $(B_1,A_1,\dots ,B_k)$ of the $kp\times kp$ matrix $J$ with $A_k>0$ such that \eqref{momrelation} holds. In this case the coefficients are necessarily unique and 
\begin{equation*}
{\mathcal{I}}_{\operatorname{m}}^{(k)}(m_1,\dots m_{2k-1}) = \mathcal{I}^{(k)}(B_1,A_1,\dots , B_{k}) . 
\end{equation*}
As in the scalar case, we do not consider the even case, since there is no injectivity there.

The Dawson-G\"{a}rtner theorem (see \cite{demboz98}) yields the LDP for the whole moment sequence $m(\Sigma^{(n)})$
in $\mathcal{H}_p^{\mathbb{N}_0}$ 
 with good rate
\begin{align*}
\mathcal{I}_{\operatorname{m}}(m_1,\dots) &= \sup_{k\geq 1} {\mathcal{I}}_{\operatorname{m}}^{(k)}(m_1,\dots m_{2k-1}) .
\end{align*}
This supremum is finite if and only if there exists a (unique) sequence $(B_1,A_1,\dots)$ of coefficients satisfying $A_k>0$ and \eqref{momrelation}. Note that this implies in particular that $(m_1,\dots)$ is the moment sequence of a nontrivial measure $\Sigma$. In this case 
\begin{align*}
\mathcal{I}_{\operatorname{m}}(m_1,\dots) &= \sup_{k\geq 1} \mathcal{I}^{(k)}(B_1,A_1,\dots , B_{k}) \\ 
& = \sup_k  \sum_{j=1}^k \frac{1}{2}\tr B_j^2  + \sum_{j=1}^{k-1}\tr G(A_j^2) \\
& = \sum_{k=1}^\infty  \left( \frac{1}{2} \tr B_k^2 + \tr\!\ G(A_k^2)\right) .
\end{align*}
$ \Box $ \\

\subsection{Laguerre block case}
\label{sLag}
The starting point for the proof of Theorem \ref{LDPhermite} is the following block-bidiagonal representation. It is a matrix extension of a famous result of Dumitriu and Edelman \cite{dumede2002}.
\begin{lem} \label{triL}
Let $m \geq n$ and for $k= 1, \dots, n$ let 
 $D_k$ and $C_k$ for $k=1,\dots ,n$ be independent random non-negative definite matrices in $\mathcal H_p$ such that  
\[C_k^2 \sim \LUE_p(p(n-k)) \ \ , \ \ D_k^2 \sim \LUE_p(p(m-k+1))\]
  and define the block matrix 
\begin{align*}
 Z_n = 
\begin{pmatrix} 
  D_1 & 0    &         &         \\
                C_1 & D_2    & \ddots  &         \\
                    & \ddots & \ddots  & 0 \\
                    &        & C_{n-1} & D_n
\end{pmatrix} .
\end{align*}
Then the $p\times p$ spectral matrix measure of $L_n = Z_nZ_n^*$ has the same distribution as the 
spectral
matrix measure of a $pn \times pn$ matrix, distributed according to the $\LUE_{pn}(pm)$ $(m\geq n)$. \end{lem}

\medskip

\proof
We use the construction of the Laguerre ensemble $L_n =   W_n W_n^*$, with $W_n$ a $pn \times pm$ matrix with independent complex Gaussian entries. Writing $w_{i,j}$ for the $p \times p$ block of $W_n$ in position $i,j$, 
let $R$ be a $pm\times pm$ unitary matrix constructed analogously to the matrix $\hat{W}\hat{H}$ in the proof of Lemma \ref{triG}, such that
\begin{align*}
W_nR  = \begin{pmatrix}
\tilde{w} \\ \widetilde W 
\end{pmatrix}
\end{align*}
with 
\begin{align*}
  \tilde w = (w_{1,1},\dots ,w_{1,m}) R = \left( \left( \sum_{i=1}^m w_{1,i}^* w_{1,i} \right)^{1/2}, 0_{p,p},\dots ,0_{p,p}\right) .
\end{align*}
The matrix $R$ can be chosen independently of $w_{i,j}, i \geq 2$ such that the entries of $\widetilde W$ are again independent complex Gaussian, independent of $\tilde w$. 
Similarly, write $z_{i,j}$ for the $p\times p$ block of $W_nR$ in position $i,j$ and let $L$ be a $ p(n-1) \times p(n-1)$ unitary matrix such that 
\begin{align*}
 L (z_{2,1}^*,\dots ,z_{n,1}^*)^* = \left( \left( \sum_{i=2}^n z_{1,i}^* z_{1,i} \right)^{1/2} , 0_{p,p},\dots ,0_{p,p}\right)^*
\end{align*} 
If $\widetilde{L} = I \oplus L$ is the extension of $L$ to an operator on $\mathbb{C}^{pn}$, leaving the subspace of the first $p$ unit vectors invariant, then
\begin{align*}
 \widetilde{L} W_nR = 
\begin{pmatrix} D_1 & 0& \dots & 0  \\
C_1 & & & \\
0   & & L\widetilde{W} R & \\
\vdots & & &  \end{pmatrix} .
\end{align*}
The first blocks satisfy $ D_1^2 \sim \LUE_p(pm), C_1^2\sim  \LUE_p(p(n-
1
))$ and by the invariance of the Gaussian distribution, the entries of $L\widetilde{W} R$ are again Gaussian distributed. 
Since we started with independent entries, all blocks in $\widetilde{L} W_nR$ are independent. The conjugation by $\widetilde{L}$ leaves the first $p$ eigenvector rows invariant, so 
$ \widetilde{L}L_n\widetilde{L}^* =  \widetilde{L} W_n R R^* W_n^*\widetilde{L}^*$ has the same spectral 
 measure as $L_n$. This yields the first step in the reduction, the iterations are straightforward.
\hfill
$ \Box $ \\

\medskip

\textbf{Proof of Theorem \ref{LDPlaguerre}:}\\
As in the proof of Theorem 
\ref{LDPhermite},
 we start by looking at the individual entries of the tridiagonal representation of Lemma \ref{triL}, now multiplied by $\frac{1}{p\gamma_n}$. For any $k$, the rescaled block $\frac{1}{p\gamma_n}C_k^2$ satisfies by Lemma \ref{basicldps} the LDP with speed $p\gamma_n$ and rate $\mathcal{I}_2$ with $\gamma=\tau$. With the speed $pn$ we would like to consider, $ \frac{1}{p\gamma_n}C_k^2$ satisfies then the LDP with rate $\tr\!\ G(\tau^{-1}\, \cdot )$ and, taking the square root, $C_k^{(n)} := \frac{1}{\sqrt{p\gamma_n}}C_k$ satisfies the LDP with speed $pn$ and rate function
\begin{align*}
\mathcal{I}_C(C) = \tr\!\ G(\tau^{-1} C^2)
\end{align*}
for $C$ positive definite and $\mathcal{I}_C(C)=\infty$ otherwise. Similarily, if we let $D_k^{(n)}:=\frac{1}{\sqrt{p\gamma_n}}D_k$, then $(D_k^{(n)})^2$ satisfies the LDP with speed $p\gamma_n$ and rate function $\mathcal{I}_2$ of Lemma \ref{basicldps} with $\gamma=1$. If we consider the speed $pn$ and the square root $D_k^{(n)}$, this changes the rate to 
\begin{align*}
\mathcal{I}_D(D) = \tau^{-1} \tr\!\ G(D^2) 
\end{align*}
for $D$ positive definite and $\mathcal{I}_D(D)=\infty$ otherwise. 

Then we follow the same way as for the Hermite model. By the independence of the matrices $C_k,D_k$, this yields the LDP for any finite sequence $(D_1^{(n)},C_1^{(n)},D_2^{(n)},\dots, D_k^{(n)}, C_k^{(n)})$ in the sequence space of Hermitian non-negative definite matrices with speed $pn$ and good rate 
\begin{equation}
\label{rateCCDD}
\mathcal{I}_{D,C} (D_1,C_1,\dots, D_k, C_k )= \sum _{j=1}^k \tau^{-1} \tr\!\ G(D_j^2) +  \tr\!\ G(\tau^{-1}C_j^2) .
\end{equation}
From (\ref{ABfromZ}), this yields the LDP for $(B_1^{(n)},A_1^{(n)},\dots ,B_k^{(n)})$ with $k$ fixed. As in the Hermite case, we may conclude a LDP for a finite collection of moments and then for the complete sequence of moments $m(\Sigma^{(n)})$ by application of the Dawson-Gärtner theorem. The resulting good rate function $\mathcal{I}_{\operatorname{m}}$ is finite in $(m_1,\dots)\in \mathcal{H}_p^{\mathbb{N}_0}$ only if $(m_1,\dots )$ is the moment sequence of a nontrivial measure with support in $[0,\infty)$. By the discussion in Section \ref{sec:polynomialson0infty}, 
this is equivalent to the existence of a sequence of positive definite matrices $D_1,C_1,D_2,\dots$ such that $(m_1,\dots)$ is the moment sequence of the spectral measure of $J=XX^*$ with $X$ as in \eqref{bidiagonal}. In this case 
\begin{equation*}
\mathcal{I}_{\operatorname{m}} (m_1,\dots)
 = \sum _{j=1}^\infty \tau^{-1} \tr\!\ G(D_j^2) +  \tr\!\ G(\tau^{-1}C_j^2) . 
\end{equation*}
We use the fact that $\tr(AB)=\tr(BA)$ and $\det(AB)=\det(BA)$ to get 
\begin{align*}
\tr\!\ G(D_k^2) = \tr\, G(D_kD_k^*) = \tr\, G(Z_{2k-1}) = \tr\, G(\zeta_{2k-1})
\end{align*}
with $Z_{2k-1}$ as in \eqref{defZZ} and \eqref{scalmat}, and 
\begin{align*}
\tr\!\ G(\tau^{-1}C_k^2) = \tr\, (\tau^{-1} C_kC_k^*)= \tr\, (\tau^{-1}Z_{2k}) = \tr\, (\tau^{-1} \zeta_{2k}) .
\end{align*}
So the value of the rate function does not depend on the unitary matrices $\sigma_n$ and $\tau_n$ in the construction of $D_k,C_k$, but only on the matrices $\zeta_k$, which in particular are uniquely determined by $(m_1,\dots )$.
The inverse contraction principle implies then the LDP for the spectral measure $\Sigma^{(n)}$. 
 \hfill
$ \Box $ \\

\section{Appendix: Extension of Baldi's theorem and Bryc's lemma}

In this part we prove a theorem which combines a LDP with a convex rate function and a LDP with a non-convex one. It is one of the  key tool for the statements in Section 4 in \cite{magicrules} and it will be used in  
 \cite{GNROPUC}.  
The first LDP deals with a random spectral measure restricted to the support of the equilibrium measure and the second LDP deals with  a subset of outliers.

To give the theorem in a general setting, assume that $\mathcal X$ and $\mathcal Y$ are Hausdorff topological vector spaces. Let $\mathcal X^*$ be the topological dual of $\mathcal X$ and equip $\mathcal{X}$
 with the weak topology. We denote by $C_b(\mathcal Y )$ the set of all bounded continuous functions $\varphi:\mathcal Y \to \mathbb{R}$. 
A point $x\in \mathcal{X}$ is called an exposed point of a function $F$ on $\mathcal{X}$, if there exists $x^*\in \mathcal{X}^*
$ (called an exposing hyperplane for $x$) such that 
\begin{align}\label{exposinghyper}
F(x) - \langle x^*,x\rangle < 
F(z) - \langle x^*,z\rangle
\end{align}
for all $z\neq x$.

\subsection{Some classical results in large deviations}

Let us recall two well known results in the theory of large deviations, which have to be combined  carefully in order to get our general theorem. The first one is the inverse of Varadhan's lemma (Theorem 4.4.2 in \cite{demboz98}), the second one is a version of the so-called Baldi's theorem  (Theorem 4.5.20 in \cite{demboz98}). The latter differs from the version in \cite{demboz98} in a straightforward condition to identify the rate function, which was applied for instance 
 in \cite{grz} (see also \cite{gamb}).
The proof of our Theorem \ref{newgeneral} will be quite similar to the 
proof of these two classical theorems. 

\begin{thm}[Bryc's Inverse Varadhan Lemma]
\label{Bryc}
Suppose that the sequence $(Y_n)$ of random variables in $\mathcal{Y}$ is exponentially tight and that the limit
\[\Lambda (\varphi) := \lim_{n\to \infty} \frac{1}{n} \log \mathbb E e^{n\varphi(Y_n)}\] exists for every $\varphi \in \mathcal C_b(\mathcal Y)$. Then $(Y_n)$ satisfies the LDP with the good rate function 
\[\mathcal I(y) = \sup_{\varphi \in \mathcal C_b(\mathcal{Y})} \{\varphi(y) - \Lambda(\varphi)\}\,.\]
Furthermore, for every $\varphi \in \mathcal C_b(\mathcal{Y})$,
\[\Lambda(\varphi) = \sup_{y\in \mathcal Y} \{\varphi(y) - \mathcal I(y)\}\,.\]
\end{thm}

\begin{thm}[A version of Baldi's Theorem]
\label{Baldi}
Suppose that the sequence $(X_n)$ of random variables in $\mathcal X$ is exponentially tight and that   
\begin{enumerate}
\item There is a  set $D\subset \mathcal X^*$ and a function $G_X:D\to \mathbb{R}$ such that for all $x^*\in D$ 
\begin{equation}
\label{ncgf}\lim_{n\to \infty} \frac{1}{n} \log \mathbb E \exp \left(n \langle x^* , X_n\rangle\right) =  G_X(x^*) \,;\end{equation}
\item The set $\mathcal F$   of exposed points $x$ of
\[G_X^* (x) = \sup_{x^* \in D} \{\langle x^*,x\rangle - G_X(x^*)\}   \]
with an exposing hyperplane $x^*$ satisfying $x^*\in D$ and 
$\gamma x^*\in D$ for some $\gamma
>1 $, is dense in 
%
$\{ G_X^* < \infty\}$.
\end{enumerate}
Then $(X_n)$ satisfies the LDP with good rate function $G_X^*$.  
\end{thm}

\subsection{A general theorem}
Our extension is the following combination of the two above theorems. The main point is that the rate function does not need to be convex, but we still only need to control linear functionals of $X_n$. 

\begin{thm}
\label{newgeneral}
Assume that $X_n \in \mathcal X$ and $Y_n \in \mathcal Y$ are defined on the same probabilistic space and 
that the two sequences $(X_n)$ and $(Y_n)$ are exponentially tight. Assume further that
\begin{enumerate}
\item There is a set $D\subset \mathcal X^*$ and functions $G_X:D\to \mathbb{R}$, $J:C_b(\mathcal Y) \to \mathbb{R}$ such that
for all $x^*\in D$ and $\varphi\in C_b(\mathcal Y)$
\begin{equation}
\label{newncgf}\lim_{n\to \infty} \frac{1}{n} \log \mathbb E \exp \left(n \langle x^* , X_n\rangle + n \varphi(Y_n)\right) =  G_X(x^*) + J (\varphi)\,;
\end{equation}
\item
The set $\mathcal F$   of exposed points $x$ of
\[G_X^* (x) = \sup_{x^* \in D} \{\langle x^*,x\rangle - G_X(x^*)\}  , \]
with an exposing hyperplane $x^*$ satisfying $x^*\in D$ and 
$\gamma x^*\in D$ for some $\gamma
>1 $, is dense in 
$\{ G_X^* < \infty\}$.
\end{enumerate}  
Then, the sequence $(X_n, Y_n)$ satisfies the LDP with speed $n$ and good rate function
\[\mathcal I(x,y) = G_X^*(x) + \mathcal{I}_Y(y)\,,\]
where 
\begin{equation*}
\mathcal{I}_Y(y) = \sup_{\varphi \in C_b(\mathcal Y)} \{ \varphi(y) - J(\varphi)\} .
\end{equation*}
\end{thm}
Let us note that in view of Varadhan's Lemma we have
\[J(\varphi) = \sup_{y \in \mathcal Y} \{ \varphi(y) - \mathcal{I}_Y(y) \}. \]

\proof

{\bf Upperbound}: The proof  follows the lines of the proof of part (b) of Theorem 4.5.3 in \cite{demboz98}. Note that since the sequence $(X_n,Y_n)$ is exponentially tight it suffices to show the upper bound for compact sets.

{\bf Lowerbound}: As usual, it is enough to consider a neighbourhood 
 $\Delta_1 \times \Delta_2$ of $(x,y)$ where  $\mathcal I(x,y) < \infty$. Take 
$\liminf_{n \to \infty} \frac{1}{n} \log \mathbb P ((X_n , Y_n) \in \Delta_1 \times \Delta_2)$ and get a lower bound tending to $\mathcal I (x,y)$ when the size of the neighbourhood tends to zero. Actually, due to the density assumption 2. it is enough to study the  lower bound of $\mathbb P(X_n \in \Delta_1, Y_n \in \Delta_2)$  when 
 $x\in \mathcal F$ and $\mathcal I_Y (y) < \infty$. 

As in \cite{demboz98} (Proof of Lemma 4.4.6), let $\varphi:\mathcal Y \rightarrow [0,1]$ be a continuous function, such that $\varphi (y) =1$ and $\varphi$ vanishes on the complement $\Delta_2^c$ of $\Delta_2$. For $m >0$, define $\varphi_m := m(\varphi - 1)$. Note that
\[J(\varphi_m) \geq -\mathcal{I}_Y(y)\,.\]

We have 
\begin{eqnarray*}
\mathbb P(X_n \in \Delta_1, Y_n \in \Delta_2) =
\mathbb E \left[ {\mathbbm{1} }_{\{ X_n \in \Delta_1\}} {\mathbbm{1}}_{\{ Y_n \in \Delta_2\}}e^{ n\langle x^*, X_n\rangle+n\varphi_m(Y_n)} e^{-n\langle x^*, X_n\rangle-n\varphi_m(Y_n)}\right] .
\end{eqnarray*}
Now $-\varphi_m \geq 0$ and on $\Delta_1$, $-\langle x^*, X_n\rangle \geq -\langle x^*, x\rangle -  \delta$ for a $\delta>0$, so that
\begin{equation}
\mathbb P(X_n \in \Delta_1, Y_n \in \Delta_2) \geq \mathbb E  \left[{\mathbbm{1} }_{\{ X_n \in \Delta_1\}} {\mathbbm{1}}_{\{ Y_n \in \Delta_2\}} e^{ n\langle x^*, X_n\rangle+n\varphi_m(Y_n)} \right] e^{- n\langle x^*, x\rangle - n \delta}\,.
\end{equation} 
Denoting
\[\ell_n = \frac{1}{n}  \log \mathbb E  e^{ n\langle x^*, X_n\rangle} \ , \ \mathcal L_n := \frac{1}{n} \log \mathbb E e^{n\langle x^*, X_n\rangle + n \varphi_m (Y_n) }\]
and 
$\widetilde{\mathbb P}$ the new probability on $\mathcal X \times \mathcal Y$ such that
\[\frac{d\widetilde{\mathbb P}}{d\mathbb P}=  e^{n\langle x^*, X_n\rangle + n \varphi_m (Y_n)  -n \mathcal L_n}\,,\]
we get
\begin{equation}
\label{mino}
\mathbb P(X_n \in \Delta_1, Y_n \in \Delta_2) \geq \widetilde{\mathbb P}(X_n \in \Delta_1, Y_n \in \Delta_2)e^{- n\langle x^*, x\rangle - n \delta + n \mathcal L_n}\,.
\end{equation}
For the exponential term we have
\begin{equation}
\label{uptodelta}
\liminf_{n\to \infty} \frac{1}{n} \log e^{- n\langle x^*, x\rangle - n \delta + n \mathcal L_n}\geq \langle x^*, x\rangle -  \delta + G_X
(x^*) + J(\varphi_m) \geq -
G_X^*(x)  -\mathcal{I}_Y(y)-\delta .
\end{equation}
We may choose $\delta $ arbitrarily small by choosing $\Delta_1$ sufficiently small, 
so that it will be enough to prove that
\begin{equation}
\label{7.7}
\widetilde{\mathbb P}(X_n \in \Delta_1, Y_n \in \Delta_2) \xrightarrow[n\to \infty]{} 1
\end{equation}
or equivalently, that
\begin{equation} \label{complementprob}
\widetilde{\mathbb P}(X_n \in \Delta_1^c) + \widetilde{\mathbb P}( Y_n \in \Delta_2^c) \xrightarrow[n\to \infty]{} 0\,.
\end{equation}
For the first term, note that under $\widetilde{\mathbb P}$ the moment generating function of $X_n$ satisfies
\begin{align*}
\lim_{n\to \infty} \frac{1}{n} \log \widetilde{\mathbb{E}} [e^{n\langle z^*,X_n\rangle}] 
& = \lim_{n\to \infty} \frac{1}{n} \log {\mathbb{E}} [e^{n\langle z^*+x^*,X_n\rangle + \varphi_m(Y_n)-n\mathcal{L}_n}] \\
& = G_X(z^*+x^*) + J(\varphi_m) - G_X(x^*) - J(\varphi_m) \\
& = G_X(z^*+x^*) - G_X(x^*)\\
& = : \widetilde G_X(z^*) ,
\end{align*}
for $z^* \in \widetilde D  := \{z^* : x^* + z^* \in D \}$. 
We may then  follow the argument on  p.159-160 in \cite{demboz98} (as an auxiliary result in  their proof of the lower bound).
Using that $x^*\in D$ is an exposing hyperplane, we get 
\begin{align*}
\limsup_{n\to\infty }\frac{1}{n}\log \widetilde{\mathbb P}(X_n \in \Delta_1^c) <0 .
\end{align*}
Considering the second term in \eqref{complementprob}, we have, on $\Delta_2^c$
\[\frac{d\widetilde{\mathbb P}}{d\mathbb P} = e^{- nm + n \langle x^*, X_n\rangle    - n \mathcal L_n }\]
so that 
\[\widetilde{\mathbb P}( Y_n \in \Delta_2^c)  \leq e^{- nm + n \ell_n    - n \mathcal L_n } \,.\]
Taking the logarithm, this implies 
\begin{align*}
\limsup_{n\to \infty} \frac{1}{n}\log \widetilde{\mathbb P}( Y_n \in \Delta_2^c) & \leq -m + 
G_X(x^*) - 
G_X(x^*)-J(\varphi_m) 
\\ & =  -m - \sup_{z\in \mathcal{Y}}\{ \varphi_m(z) - \mathcal{I}_Y(z)\} \leq -m + \mathcal{I}_Y(y)
\end{align*}
which tends to $-\infty$ when $m \rightarrow \infty$.

To summarize, we have proved (\ref{complementprob}), i.e. (\ref{7.7}), which with (\ref{mino}) and (\ref{uptodelta}) gives
\[\lim_{\Delta_1 \downarrow x, \Delta_2 \downarrow y}\liminf_{n \to \infty} \frac{1}{n}\log \mathbb P(X_n \in \Delta_1, Y_n \in \Delta_2) \geq - G_X^*(x) - \mathcal I_Y (y)\,,\]
which leads to the lower bound of the LDP.
\hfill $\Box$

\begin{rem}
In Section 4 of \cite{magicrules}, and similar to Section 5.1 of the present paper, we studied the joint moment generating function of 
$(\tilde\mu^{(n)}_{I}, \lambda_M^\pm )$. 
For $s \in \mathbb R^{2j}$ we introduced
\begin{align*}
 G_n (f,s) = \mathbb E \left[\exp \left\{n \int f \, d\tilde\mu^{(n)}_{I} + n \langle s,\lambda_M^\pm\rangle \right\}\right] ,
\end{align*}
and proved that 
 for all $f$ such that $\log (1-f)$ is continuous and bounded (and all $s\in \mathbb{R}^{2j}$), 
\begin{align} \label{oldlimit}
\lim_{n\to \infty} \frac{1}{n}\log G_n (f,s) = G(f) + H(s)\,.
\end{align}
Actually $H$ is the Legendre dual of $\mathcal I_{M, \lambda^\pm}$, i.e.
\[H(s) = \sup_{y \in \mathbb R^{2j}} \{ \langle s, y\rangle - \mathcal{I}_{M,\lambda^\pm}(y)\} = \left(\mathcal{I}_{M,\lambda^\pm}\right)^*(s) \,.\]


However, the rate $\mathcal{I}_{M,\lambda^\pm}$ might be non-convex (when $V$ is not convex) and hence the dual $H^*$ is not strictly convex on a dense set, and then the Assumption 2 of Theorem \ref{Baldi} may not be verified. The convergence in \eqref{oldlimit} is therefore not enough to conclude the joint LDP for $(\tilde\mu^{(n)}_{I},\lambda_M^\pm)$ directly from the  Theorem \ref{Baldi} above. A complete proof of Theorem 3.1 in \cite{magicrules} needs actually an application of Theorem \ref{newgeneral}.

\end{rem}
\bibliographystyle{alpha}
\bibliography{bibclean}

\end{document}